\documentclass[a4paper, 11pt]{article}
\usepackage{amssymb,latexsym,amsmath,amsthm,amsfonts}
\usepackage{tikz}
\usepackage{paralist}	 
\usepackage{cite}
\usepackage{color}
\usepackage{multicol}
\usepackage{caption}
\usepackage{subfigure}
\usepackage[latin1]{inputenc}
\usepackage{graphicx}
\usepackage{epsf}
\usepackage{fancyhdr}
\usepackage{pstricks}
\usepackage{algorithm}
\usepackage{algorithmic}
\usepackage{lscape}
\newpsobject{showgrid}{psgrid}{subgriddiv=1, griddots=10,gridlabels=6pt}
\usepackage{authblk}
\usepackage[capitalise]{cleveref}

\title{Packing a fixed number of identical circles in a circular container with circular prohibited  areas}   
\author[1]{C.O. López\thanks{claudia.lopez@ciencias.unam.mx}}
\author[2]{J.E. Beasley\thanks{john.beasley@brunel.ac.uk; john.beasley@jbconsultants.biz}}

\affil[1]{Faculty of Sciences, National Autonomous University of Mexico, Mexico City, Mexico}
\affil[2]{Mathematical Sciences, Brunel University, Uxbridge UB8 3PH, UK and JB Consultants, Morden, UK }

\begin{document}
\date{June 2017, Revised July 2018}   
\par\medskip

\maketitle

\begin{abstract}
In this paper we consider the problem of packing a fixed number of identical circles inside the unit circle container, where the packing is complicated by the presence of fixed size circular prohibited areas. Here the objective is to 
maximise the radius of the identical circles.
We present a heuristic for the problem based upon formulation space search. 
Computational results are given for  six test problems involving the packing of up to 100 circles. One test problem has a single prohibited area made up from the union of circles of different sizes. Four test problems are  annular containers, which have a  single inner circular prohibited area. One test problem has circular
prohibited areas that are disconnected.
\end{abstract}


\section{Introduction} \label{Sec:Intro} 
In this paper we consider the problem of packing a  given number of circles of identical size inside a circular container. This problem can be viewed from two different, but equivalent, perspectives: 
\begin{enumerate}
\item minimise the size of the container having fixed size circles to be packed; or
\item maximise the identical size of the circles to be packed having a container of fixed size.
\end{enumerate}
These problems are equivalent since there is a simple relationship between them based on scaling. 

\begin{figure}[!htbp] 
 \centering 
\subfigure[Prohibited area: Union of circles]{\label{fig:fa1}\includegraphics[scale=0.3]{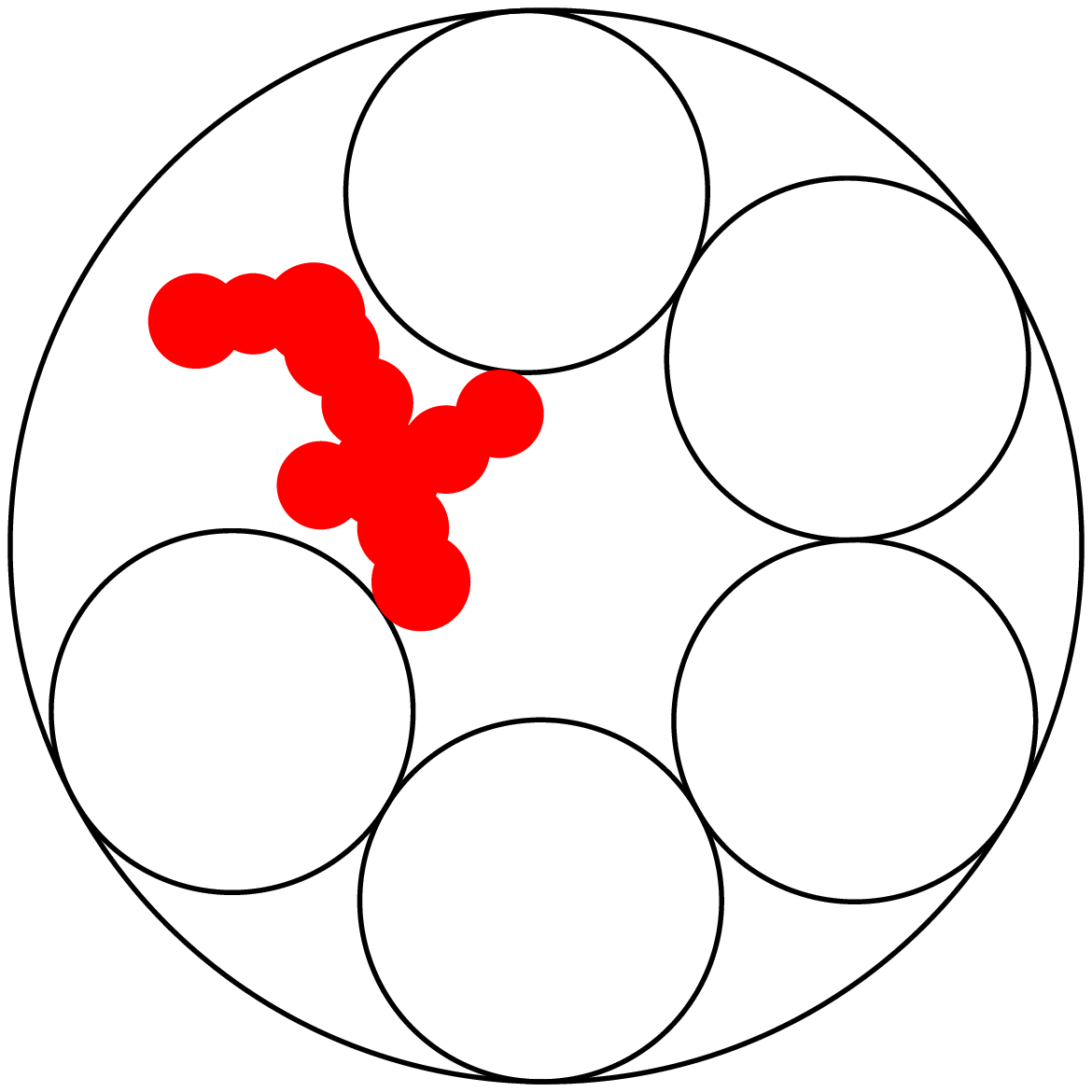}} \hspace{.5cm}
\subfigure[Prohibited area: Centre placed annular]{\label{fig:fa2}\includegraphics[scale=0.3]{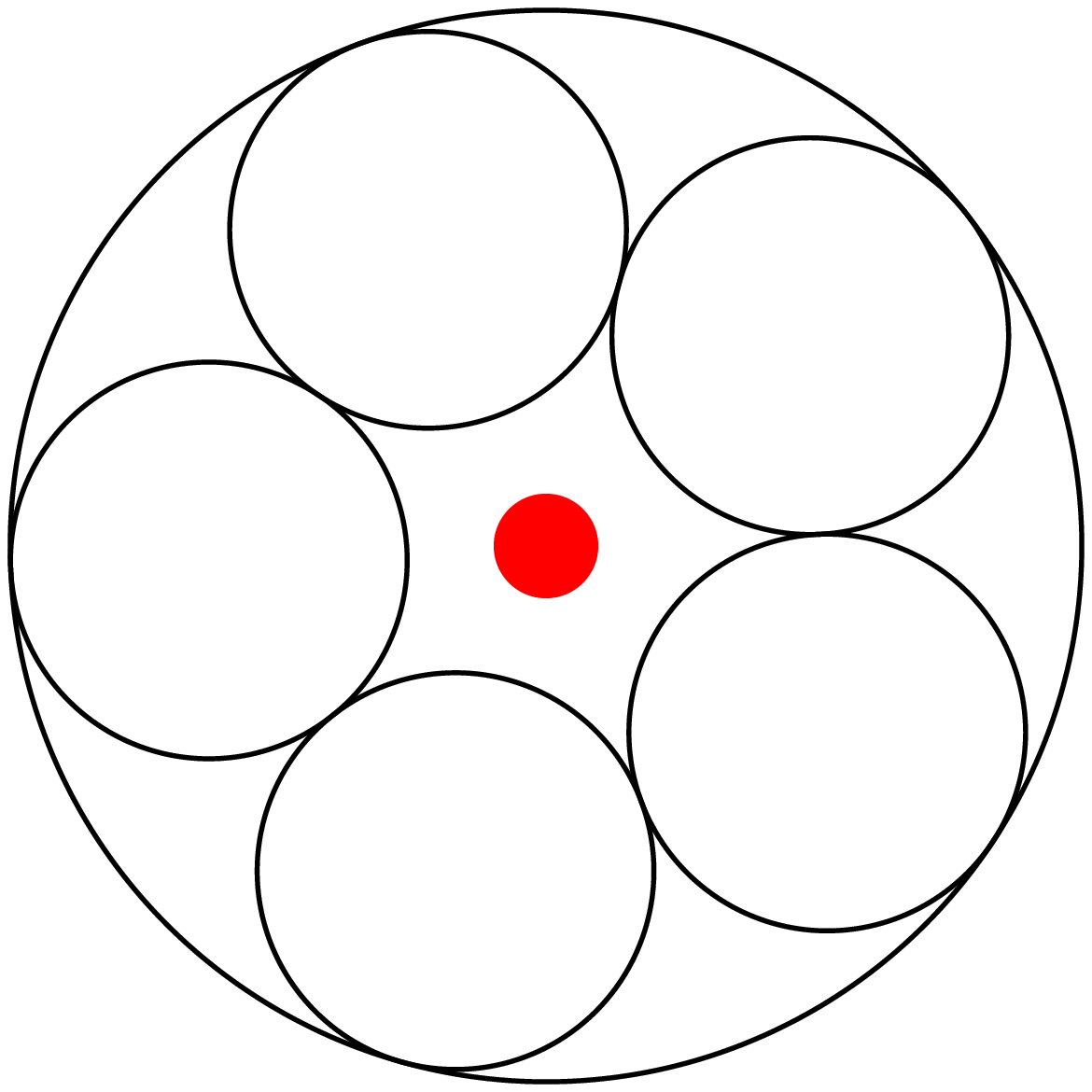}}\hspace{.5cm}
\subfigure[Prohibited area: Bottom placed annular]{\label{fig:fa3}\includegraphics[scale=0.3]{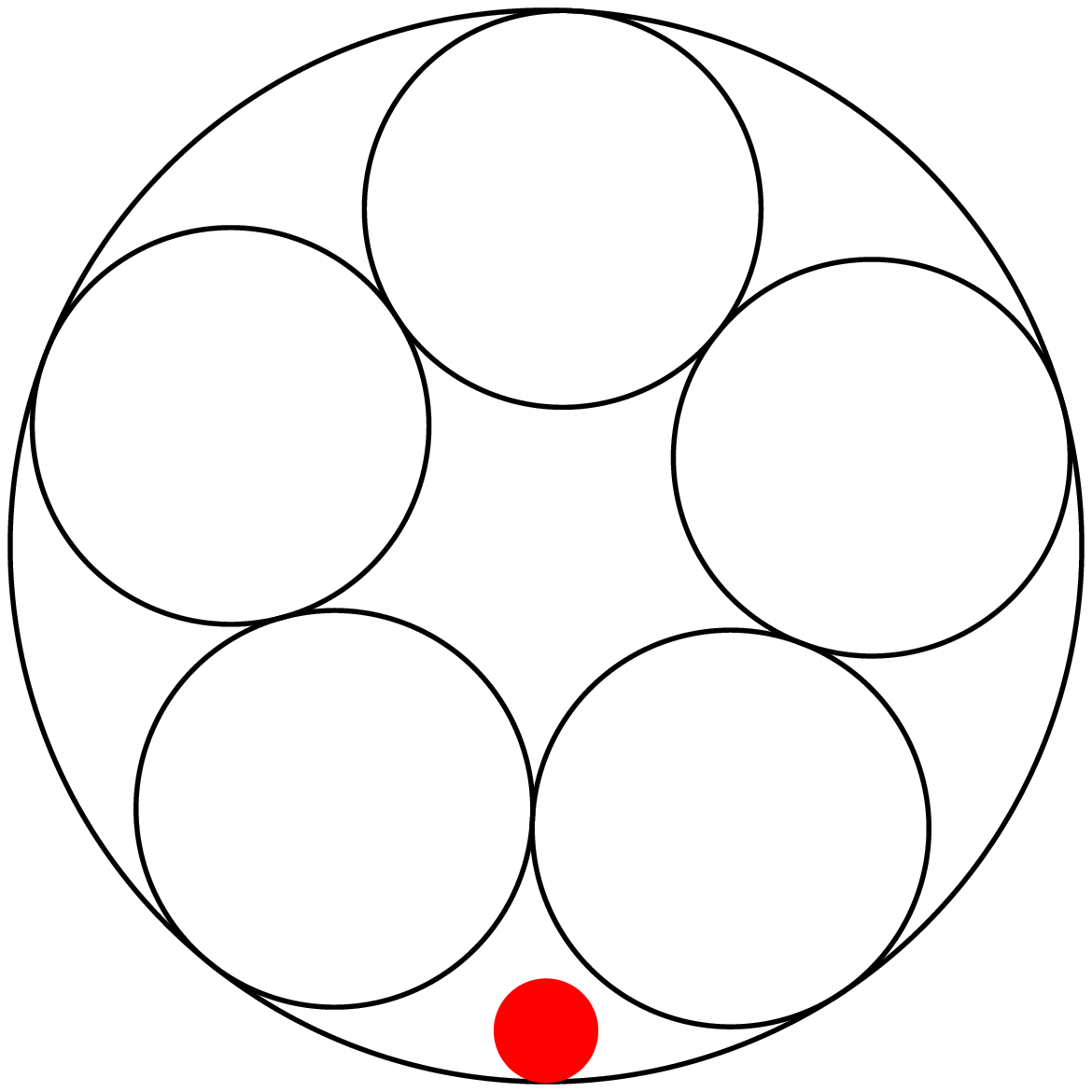}}
\caption{Circle packing problem with $n=5$ identical circles with circular prohibited areas}
\label{Fig:5circles} 
\end{figure}

In this paper we adopt the second perspective,  with the container being the unit circle. We  consider  problems with differing types of fixed size circular prohibited areas. The first type has a single prohibited area as the union of circles of different  sizes, as shown 
in~\cref{fig:fa1}. The second type  are  annular containers, which have a  single inner circular prohibited area, such 
 as those shown in~\cref{fig:fa2,fig:fa3}. In~\cref{fig:fa1,fig:fa2,fig:fa3} we show a packing of five identical circles, where the prohibited areas are 
shown as solid circles. 
The third type has circular prohibited areas that are disconnected.

Note here that a characteristic of the packing problem considered in this paper is that, as the number of circles to be packed is fixed, and the objective is to maximise the (identical) radius of these circles, a packing might well involve unused space. This behaviour can be seen in 
in~\cref{fig:fa1,fig:fa2,fig:fa3}. In particular note that the packings seen in~\cref{fig:fa2,fig:fa3}, which are known to be optimal from geometric considerations, have a significant amount of unused space. Indeed none of the circles seen 
in~\cref{fig:fa2,fig:fa3}
are even in contact with the prohibited area.

The heuristic presented here is an iterative procedure that uses  formulation space search (henceforth FSS). This method can be seen as a  procedure in which different formulations of the problem are used at every iteration in order to try and discover better solutions. The way we construct the formulations is explained below. The contribution of this paper is:
\begin{compactitem} 
\item to consider the problem of maximising the radius of a given number of identical circles when packed inside a circular container with fixed size circular prohibited areas
\item to present computational results for  six test problems, involving differing prohibited areas,
 with up to 100 circles
\item to use formulation space search, a new and emerging metaheuristic
\end{compactitem} 

\noindent Although the prohibited areas seen in 
in~\cref{fig:fa1,fig:fa2,fig:fa3} are connected (albeit trivially so in the case of~\cref{fig:fa2,fig:fa3}) this is not a restriction on the heuristic presented in this paper. We do present below computational results for a test problem where the prohibited  areas are disconnected.

The motivation for studying the circle packing problem with prohibited areas considered in this paper is that the problem  appears when adopting a sequential/priority approach to packing. In some situations certain circular items have to be packed first, followed by a packing of the remaining circular items. In terms of the problem considered in this paper the  circular items that are packed first become prohibited areas with regard to the packing of the remaining circular items.

This paper is organised  as follows. In Section~\ref{Sec:LitSurvey} we present a literature survey, reviewing  circle packing with prohibited areas and its practical applications and FSS. In Section~\ref{Sec:Formulation} we present our formulation of the problem whilst 
Section~\ref{Sec:Heuristic} describes the algorithm proposed. 
In Section~\ref{Sec:Comparisons} we present  computational results  and finally in Section~\ref{Sec:Conclusions} we present conclusions.
 
\section{Literature survey}\label{Sec:LitSurvey}
In this section we review the work that has been done with respect to circle packing and prohibited areas. We discuss the practical applications of the circle packing problem with prohibited areas considered in this paper. 
We also give a brief introduction to FSS and the work that has been done using FSS. Note here that circle packing without prohibited areas has been extensively discussed in the literature (e.g.~see the survey papers~\cite{Castillo2008, Hifi2009}) and so, for reasons of space, we do not discuss that work here.

\subsection{Circle packing with prohibited areas}

The circle packing problem with circular prohibited areas has not been considered at all in the literature until comparatively 
recently.

Stoyan and Yaskov~\cite{Stoyan2012} considered the problem of packing the maximal number of identical circles of fixed radius inside a fixed size circular container with fixed size circular prohibited areas. Their approach, based on work in~\cite{Stoyan2012a}, defines a function for each circle $i$ that is one if circle $i$ is inside the container and does not intersect with any prohibited area, zero otherwise. They then maximise the sum of these functions subject to the constraint that the circles do not intersect. Their objective function is therefore discontinuous and takes values $0,1,2, \ldots, n$.  This is then reduced to a finite sequence of problems with linear objective functions.
In order to generate a feasible initial solution they 
use a modification of the Zountendijk method of feasible directions. Detailed computational results were given for 13 test problems.

In~\cite{Stoyan2012} Stoyan and Yaskov 
refer to the earlier work 
of Li and Akeb~\cite{Li2005}. Although~\cite{Li2005}  does not refer to prohibited areas Li and Akeb do consider the problem of packing the maximal number of identical circles of fixed radius inside a fixed size circular container. In particular, they considered an annular container, where there is an single inner fixed circle that cannot be used, so the packed circles cannot intersect with this inner circle. In other words the inner circle is a single prohibited area, as shown for example in~\cref{fig:fa2,fig:fa3}.
Li and Akeb~\cite{Li2005} presented two greedy heuristics based upon differing rules for circle placement. With  respect to annular containers detailed computational results were given for 9 test problems.
Note here that the problem considered in ~\cite{Li2005,Stoyan2012} is \textbf{\emph{different}} from the problem considered in this paper. In~\cite{Li2005,Stoyan2012}  the circle radius is fixed and the optimisation relates to the number of circles packed. In our work in this paper the number of circles is fixed and the optimisation relates to the circle radius.

With respect to packing circles with prohibited areas the only other work we are aware of is~\cite{Stoyan2012b} which deals with minimising the length of  a strip (which contains prohibited areas) when packing both circles and non-convex polygons. 

With regard to a related problem involving prohibited areas 
Zhuang et al~\cite{Zhuang2015} considered the problem of packing a fixed number of unit circles inside a square so as to minimise the size of the containing square (minimising the length of the side of the square, thereby also minimising its area and circumference). In their approach the problem is complicated by the presence within the square of a number of prohibited areas (which they term damaged areas) that cannot overlap with any unit circle. In the computational examples they considered these damaged areas were composed of very small non-overlapping squares that, in total, only constituted approximately 2.22\% of the area of the containing square. They presented an approach based on an enhancement of greedy vacancy search~\cite{Huang2010} and too used simulated annealing. Computational results were presented for a number of randomly generated problems involving the packing of up to 70 unit circles.

\subsection{Applications}
The circle packing problem has a long history and a wide variety of applications. An introduction to 
its history can be found in
Szabó et al~\cite{Szabo2007}. 
Industrial applications of the circle packing problem 
such as: circular cutting, container loading, cylinder packing, 
facility dispersion and communication networks, facility and dashboard layout,  
are considered in Castillo et al~\cite{Castillo2008}. Other applications that have been considered in the literature include the 
packing of optical fibres
into tubes, see Wang et al~\cite{Wang2002}.

With regard to the circle packing problem with prohibited areas considered in this paper then, as mentioned previously above, the problem  appears when adopting a sequential/priority approach to packing. So here, for example for a number of practical reasons, certain circular items have to be packed first, followed by a packing of the remaining circular items. In terms of the problem considered in this paper the  circular items that are packed first become prohibited areas with regard to the packing of the remaining circular items.

For example, consider~\cref{fig:fa3} and regard the circular container as a tube into which we have to pack a number of circular cables.  For a variety of reasons one particular cable has to be placed at the bottom of the tube. One reason why this might occur is simply weight considerations, if this cable is much heavier than the other cables to be positioned in the tube it naturally makes sense to have it at the bottom of the tube. So this cable has to be positioned first, and hence becomes a prohibited area for the packing of the remaining circular cables (i.e.~we have a bottom placed annular prohibited area as seen in~\cref{fig:fa3}).
We would note here that the algorithm presented in this paper is capable of dealing with any number of prohibited areas, wherever they are positioned in the circular container.

Another practical example where some  circular items have to be packed first, followed by a packing of the remaining circular items, would relate to situations (as in Wieman~\cite{Wieman2008}) where certain of the circular items are tubes 
through which coolant flows
and have to be prepositioned in order to provide adequate coolant coverage for the other items that have yet to be packed.

A recent practical application that adopts a sequential packing approach  is given in Pedroso et al~\cite{Pedroso2014}. In that paper the authors consider the  problem of the packing of tubes within a single circular outer tubular container. 

\subsection{Formulation space search (FSS)}
When solving nonlinear non-convex problems with the aid of a solver, 
Mladenovi\'{c} et al~\cite{Mladenovic2005} observed that different formulations of 
the same problem may have different characteristics.  
Hence a natural way to proceed is by
swapping between formulations. Under this framework Mladenovi\'{c} 
et al~\cite{Mladenovic2005} used FSS for the circle packing problem considering two formulations of the problem: one in a 
Cartesian coordinate system, the other in a Polar coordinate system. Their algorithm solves the problem with one formulation 
at a time and when the solution value is the same for both formulations  the algorithm terminates. They considered the case 
of packing identical circles and their  computational results were for up to 100 circles packed into the unit circle and the unit square. 

In Mladenovi\'{c} et al~\cite{Mladenovic2007} they improved 
on~\cite{Mladenovic2005} by considering a mixed 
formulation of the problem; they set a subset of the circles in the Cartesian system whilst the rest of 
the circles were in the Polar system. A reduction in the number of the non-overlapping constraints was made 
at the initial solution by disregarding points sufficiently far away from each other. They gave computational 
results for up to 100 identical circles inside the unit circle.
López and Beasley~\cite{BL2011} used FSS for the problem of packing 
equally sized circles inside a variety of containers. They presented computational results which 
show that their approach improves upon previous results based on FSS presented in the literature.
For some of the containers considered they improve on the best result known previously.
López and Beasley~\cite{BL2013} used FSS  to solve the packing problem with non-identical circles in different shaped containers. They presented computational results which were compared with benchmark problems and also proposed some new instances.

López and Beasley~\cite{BL2016} used FSS  to solve the problem of packing  non-identical circles in a fixed size circular container. This involves making a choice as to which circles to pack, as well as deciding the location of any packed circle within the container.
López and Beasley~\cite{BL2018} used FSS  to solve the problem of packing  unequal rectangles and squares in a fixed size circular container. This also involves making a choice as to which objects to pack, as well as deciding the location of any packed object
 within the container. With respect to the packing of rectangles they considered both 
fixed orientation and rotation through   90 degrees.

Essentially underlying FSS is the fact that because of the nature of the solution process in nonlinear optimisation we often fail to obtain a globally optimum solution from a single formulation. This may be because of the non-convex nature of the formulation considered, or due to user termination of the search process for computational reasons (e.g.~terminating the search upon reaching a predefined time limit).  As a consequence of this 
perturbing/changing the formulation and resolving the nonlinear program may lead to an improved solution. This leads naturally to the idea
of constructing iterative schemes that move between formulations in a systematic manner.

Note here that within FSS we consider a number of different formulations of the same problem. However FSS does not need to define a distance metric (distance function) that measures the distance between any two formulations. Note too here that FSS is distinctly different from repeatedly resolving a single formulation with different random starting solutions. Rather FSS, by considering different formulations of the same problem, potentially provides more opportunities to obtain good solutions.

FSS has been applied to a few problems  additional to circle packing (e.g.~timetabling~\cite{ Kochetov2008}). In~\cite{BL2014}  FSS was used   to solve some benchmark mixed-integer nonlinear programming problems. In a more general sense an adaptation to  FSS was presented in~\cite{Brimberg2014} for solving continuous location problems. More discussion as to FSS can be found in Hansen et al~\cite{Hansen2010}. A related approach is variable space search, which has been applied to graph colouring (Hertz et al~\cite{Hertz2008, Hertz2009}). Other related approaches are variable formulation search which has been applied to the cutwidth minimisation 
problem~\cite{Pardo2013,Duarte2016}
and variable objective search which has been applied to  the maximum independent set problem~\cite{butenko13}.

As noted in Pardo et al~\cite{Pardo2013}  
variable space search, 
variable formulation search and 
variable objective search 
contain similar ideas as originally expounded using FSS. At a slightly more general level FSS can be regarded as a variant of variable neighbourhood search, for example 
see~\cite{Amirgaliyeva2017, Hansen2017}.
Recently Erromdhani et al~\cite{err2017} used the phrase variable neighbourhood formulation search to describe their approach for solving the multi-item capacitated lot-sizing problem with time windows and setup times. In their approach both formulations and neighbourhoods are changed during the search process.

\section{Formulation} \label{Sec:Formulation}
The formulation that we use to address the packing of $n$ circles with identical size into a two-dimensional circular container (of unit radius, centred at the origin) with fixed size circular prohibited areas
involves the following notation. Let:
\begin{compactitem} 
\item $C$ be the (indexed) set of circles whose centres are expressed in 
Cartesian coordinates, so for circle $i \in C$ its centre is at $(x_i,y_i)$ in Cartesian coordinates
\item $P$ be the (indexed) set of circles whose centres are expressed in Polar coordinates, 
so for circle $i \in P$ its centre is at $(r_i,\theta_i)$ in Polar coordinates, 
where $P\cap C =\emptyset$ and $P\cup C = \{1,...,n\}$
\item $Q$ be the set of all pairs $\{(i,j)| i=1,...,n; j=1,...,n; i\ne j\}$
\item $F$ be the set of fixed prohibited circular areas, where circle $f \in F$ has radius $R_f^*$, being centred at $(x_f^*,y_f^*)$ in Cartesian coordinates (correspondingly $(r_f^*,\theta_f^*)$ in Polar coordinates)
\item $Q^*$ be the set of all pairs $\{(i,f)| i=1,...,n; f \in F \}$
\item $R$ be the radius associated with each of the $n$ circles 
\item $R_{overall}$ be any suitable upper bound on $R$, 
defined here 
by equating the area of the containing unit circle minus the area of the largest circular prohibited area ($\pi 1^2 - max_{f\in F} \pi (R_f^*)^2$) to the total area of the $n$
circles ($n \pi R_{overall}^2$), so $R_{overall}=\sqrt{(1-max_{f\in F}(R_f^*)^2)/n}$
\end{compactitem} 

\noindent Although we have above (using disjoint sets $C$ and $P$) separated centres expressed in Cartesian and Polar coordinates, note here that the relationship between the two coordinate systems is that a point $(x,y)$ in Cartesian space has equivalent coordinates 
$(r,\theta)$ in Polar space where $x=r\cos(\theta)$ and 
$y=r\sin(\theta)$.

 Our formulation of the problem then is:
\vspace{-0.5cm}
\allowdisplaybreaks 
\begin{center}
\begin{align} 
\max  & \hspace{.5cm} R & & \label{e1}\\ \notag
\textrm{subject~to} &  & &\\  
  &x_i^2 + y_i^2 \leq  (1-R)^2, &  \forall &i\in C \label{e2}\\
  &r_i \leq  1-R, &  \forall &i\in P \label{e3}\\
  &(x_i-x_j)^2 + (y_i-y_j)^2 \geq  4R^2, &  \forall (i,j)\in Q \hspace{.1cm}\text{with} &\hspace{.1cm} i\in C \hspace{.1cm} j\in C \hspace{.1cm}i<j \label{e4}\\
  &(x_i- r_j\cos(\theta_j))^2 + (y_i- r_j\sin(\theta_j))^2 \geq 4R^2, &  \forall (i,j)\in Q \hspace{.1cm} \text{with} &\hspace{.1cm} i\in C \hspace{.1cm} j\in P \label{e5}\\
  &r_i^2 + r_j^2 -  2r_i r_j \cos(\theta_i - \theta_j) \geq  4R^2, & \forall (i,j)\in Q \hspace{.1cm} \text{with} &\hspace{.1cm} i\in P \hspace{.1cm} j\in P \hspace{.1cm} i<j \label{e6}\\  	
	&(x_i-x_f^*)^2 + (y_i-y_f^*)^2 \geq  (R + R_f^*)^2, & \forall (i,f)\in Q^* \hspace{.1cm} \text{with} &\hspace{.1cm} i\in C \hspace{.1cm} f\in F \hspace{.1cm} \label{e7}\\
	&r_i^2 + (r_f^*)^2 - 2r_i r_f^* \cos(\theta_i - \theta_f^*)\geq (R + R_f^*)^2, &   \forall (i,f)\in Q^* \hspace{.1cm} \text{with} &\hspace{.1cm} i\in P \hspace{.1cm} f\in F \hspace{.1cm} \label{e8}\\	
  &-1 \leq x_i \leq 1,  &\forall &i\in C \label{e9}\\
  &-1 \leq y_i \leq 1,  & \forall &i\in C \label{e10}\\
  &0 \leq r_i \leq 1,  &\forall &i\in P \label{e11}\\
  &0 \leq \theta_i \leq 2\pi,  & \forall &i\in P \label{e12}\\  
  &0 \leq R \leq R_{overall}   & & \label{e13}
\end{align} 
\end{center}

\noindent The objective function, Equation~(\ref{e1}), maximises the radius of the identical circles to be packed. Constraints~(\ref{e2}) and~(\ref{e3}) ensure that the centre coordinates of any circle lie inside the container, in this case the unit circle. Constraint~(\ref{e4}) imposes the non-overlapping condition, that is that for any pair of circles with centres coordinates $(x_i,y_i)$ and $(x_j,y_j)$ the distance between them is at least  $2R$. In a similar manner constraints~(\ref{e5})-(\ref{e6}) are the non-overlapping condition if one or both circles are in the Polar coordinate system. Constraints~(\ref{e7})-(\ref{e8}) guarantee that no  circle will overlap any circle in the  prohibited area.  Constraints~(\ref{e9})-(\ref{e12}) represent the bounds for the variables in the Cartesian and Polar coordinate system, while constraint~(\ref{e13}) imposes bounds on the variable $R$.
We would stress here that in our formulation the prohibited areas need not be connected together, and there can be as many of them as we wish.

In the formulation presented above it is clear that 
Equations~(\ref{e9}) and~(\ref{e10}), which impose bounds on the Cartesian coordinates, can be deduced from Equation~(\ref{e2}). We have presented these bounds here however for completeness.  A similar remark applies to the upper limit on the Polar radius in 
Equation~(\ref{e11}), which can be deduced from Equation~(\ref{e3}).

We would mention here that there is an alternative approach to specifying the non-overlapping conditions (Equations~(\ref{e4})-(\ref{e8})) seen above based on that presented in Birgin and  Sobral~\cite{Birgin2008}. This involves replacing the (numerous) individual non-overlapping constraints by a single constraint involving the summation of maximisation terms. We did investigate that approach, but were unable to successfully implement it. This was due to software limitations, in particular with regard to incorporating the summation of maximisation terms within a constraint in the nonlinear solver  (SNOPT~\cite{snopt,snoptTom}) which we used.

\section{Heuristic} \label{Sec:Heuristic}
The heuristic developed is an iterative procedure that uses FSS. At each iteration we first define a new formulation of the problem; then solve it using a nonlinear solver  (SNOPT~\cite{snopt,snoptTom}); then apply a correction step to account for possible numerical inaccuracies. 
Our FSS heuristic is an extension of that given in~\cite{BL2011}, which itself built upon the earlier work of~\cite{ Mladenovic2005, Mladenovic2007}.

\subsection{Formulation}\label{section:swap}
Suppose that we have current Cartesian coordinate positions for the $n$ circles given by 
$(X_i,Y_i)~i=1,\ldots,n$. We first restrict the problem such that circles $i \in C$ cannot be positioned more than a user defined distance $\Delta$ from their current position in a horizontal or vertical direction when we come to resolve the problem. To achieve this we add to the formulation given previously the constraints:
\allowdisplaybreaks 
\vspace{-0.5cm}
\begin{center}
\begin{align} 
  &X_i-\Delta \leq x_i \leq X_i+\Delta,  &\forall &i\in C \label{e14}\\
  &Y_i-\Delta \leq y_i \leq Y_i+\Delta,  & \forall &i\in C \label{e15}
\end{align} 
\end{center}

Stoyan et al~\cite{Stoyan2016} use similar constraints to  
Equations~(\ref{e14})-(\ref{e15}), but in the context of packing ellipses. 
We have also used similar constraints in~\cite{thesis, BL2011}. 
Based on work presented in~\cite{thesis, BL2011} we only restrict the positions of circles whose centres are expressed in Cartesian coordinates.  This allows additional flexibility for the positioning of circles
whose centres are expressed in Polar coordinates.

Now given the above restriction on the movement of circle centres for circles $i \in C$ it is  a simple matter to update the set $Q$ of circle pairs to eliminate from that set any pairs of circles $(i,j)$ which, because of the restriction on the movement of the centres of $i \in C$ and $j \in C$ (and the restriction $R_{overall}$ on maximum circle size), can never overlap. This enables us to have fewer constraints when we come to solve our formulation. In the same manner we can also update the set $Q^*$ of circle pairs to eliminate from that set any pairs of circles $(i,f)$ which, because of the restriction on the movement of the centre of $i \in C$, can never overlap. 

As will become apparent in the pseudocode given below, we also at each iteration randomly allocate circles $i,~i=1,\ldots, n$ to $C$ or $P$ (by for each circle $i$ in turn taking a random number from the uniform distribution $[0,1]$ and assigning circle $i$ to $C$ if the random number is 0.5 or less, assigning it to $P$ otherwise).

At each iteration therefore this perturbation of the circles assigned to $C$ and $P$, together with the additional constraints 
(Equations~(\ref{e14})-(\ref{e15})) on movement, means that we have a \textbf{\emph{different nonlinear formulation}} of the problem to solve.  

We commented before that FSS is distinctly different from repeatedly resolving a single formulation with different random starting solutions. Rather FSS, by considering different formulations of the same problem, here for example as we change the sets $C$ and $P$ at each iteration, potentially provides more opportunities to obtain good solutions.

Note here that our formulation will give a heuristic result because of the constraints on the movement of circle centres. But the logic here is that we cannot realistically expect to solve the original problem to global optimality anyway, so imposing such constraints 
(reducing the size of the problem, here by reducing the size of both $Q$ and $Q^*$, thereby reducing the number of constraints that need be considered, Equations~(\ref{e4})-(\ref{e8})) 
may have computational benefits. In fact, as discussed later below, our computational experience has been that introducing 
Equations~(\ref{e14})-(\ref{e15}) does significantly improve the quality of the results that we obtain.

We solve our formulation of the problem as presented in Equations~(\ref{e1})-(\ref{e15}), denoted as  $NLP(C,P)$ in the pseudocode given below, using the  nonlinear solver  
SNOPT~\cite{snopt,snoptTom}.

\subsection{Correction procedure}~\label{Subsec:Correction}
Even though the result given by the nonlinear solver  has a high degree of accuracy we have included a
correction procedure to avoid numerical inaccuracies. 
In the literature (e.g.~see~\cite{Pack}) solutions
to circle packing problems are typically given to a high degree of numeric accuracy and as such it is appropriate to regard the circle centre coordinates as given by the nonlinear solver as fixed and to adjust the circle radius $R$ to ensure that all constraints are satisfied to a high degree of numeric accuracy. Simply taking the value for $R$ as output by the solver may not be sufficient. Issues related to numeric accuracy (in the absence of a correction procedure such as that outlined below) become particularly acute as $n$ increases and the circles are packed closer and closer together.

As we are solving the circle packing problem with prohibited areas all solutions must satisfy three conditions: no circles overlap with a prohibited area, all circles are inside the circular container and no overlapping circles. Setting a value for $R$ to ensure that these three conditions are met is easily done and so, for space reasons, we omit details here. Examples of correction procedures can be found in our earlier work~\cite{thesis, BL2011, BL2013}.

We would stress here that applying a correction procedure 
is especially important in circle packing since there are examples in the 
literature, e.g.~\cite{Birgin2008} as noted in~\cite{Birgin2010}, where the results given are 
invalid due to a lack of sufficient precision. In the computational implementation of our correction procedure we used a 
MATLAB function called vpa (which is the acronym for variable precision arithmetic) 
that gives as many digits of accuracy as we desire. Note though that we use default machine precision in all other computations.

\subsection{Pseudocode}
The pseudocode is shown in Algorithm~\ref{colsalg}. The first initial solution is generated randomly inside the container, in our case the unit circle, by choosing $r_i$ uniformly from the interval $[0,1]$ and $\theta_i$ uniformly from interval $[0,2\pi]$ and expressing their respective values in Cartesian coordinates $(X_i,Y_i)$. 
At each iteration after updating the sets $Q$ and $Q^*$ we solve the nonlinear program (Equations~(\ref{e1})-(\ref{e15})) using  
SNOPT~\cite{snopt,snoptTom}. We correct the resulting radius to $R^*$ and update the best radius found ($R_{best}$) accordingly. 
We also update the value for $\Delta$, where we  set it to 
$\frac{2}{3}R^*$, so relate it to the value of the current solution. This value of $\frac{2}{3}$ was based on limited computational experimentation carried out as reported in~\cite{thesis}.
We then update the iteration counter, the current solution and the circle sets 
$C$ and $P$ and repeat until termination.  Note here that we update the current solution  even if the solution has not improved as compared with the preceding solution.

In the computational results reported below the termination criteria was set to 80 iterations and we performed 25 replications, reporting the best (maximum) radius found over all replications. This means that for each value of $n$ examined we solve $80(25)=2000$ nonlinear programs.

\begin{algorithm}[!htb]
\phantom{a}
\caption{Pseudocode for the heuristic}
\begin{algorithmic}
\STATE \textbf{Initialisation:} 
\STATE $t \leftarrow 0$;\ $R_{best} \leftarrow 0$; \ $\Delta \leftarrow \frac{2}{3}R_{overall}$ \\ 
\STATE Randomly generate an initial solution $(X_i,Y_i)~i=1,\ldots,n$  \\
\WHILE{not termination condition} 
\STATE Update $Q$ and $Q^*$  \hfill  \COMMENT{Update the sets of pairs}\\   
\STATE $(x,y,R) \leftarrow$ Solve $NLP(C,P)$
\hfill \COMMENT{Solve the formulation} \\
\STATE  $R^* \leftarrow$  value after correction \hfill \COMMENT{Correct the radius} \\
\STATE $R_{best} \leftarrow \max \{R_{best},R^*\}$ \hfill \COMMENT{Update the best radius $R_{best}$}  \\
\STATE $\Delta \leftarrow \frac{2}{3}R^*$  \hfill \COMMENT{Update $\Delta$}\\
\STATE $t \leftarrow t+1$ \hfill \COMMENT{Update the iteration counter}\\
\STATE $(X,Y) \leftarrow (x,y)$ \  \hfill \COMMENT{Update the current solution}\\
\STATE Update $C$ and $P$   \  \hfill \COMMENT{Update the circle sets by randomly allocating circles to $C$ or $P$}\\
\ENDWHILE
\end{algorithmic}
  \label{colsalg}
\medskip
\end{algorithm}

\section{Computational results} \label{Sec:Comparisons}
The results presented in this section for our FSS heuristic were produced on an Intel(R) Core(TM) i5-2500 3.30GHz CPU with 4GB of memory. The algorithm was coded in MATLAB 7.9.0 
using SNOPT~\cite{snopt,snoptTom} as the nonlinear solver.
We considered six different test problems, where the details of the prohibited areas associated with each test problem can be found in Table~\ref{Table:coordinates}.

In  Table~\ref{Table:coordinates} test problem 1 
has a single prohibited area made up from the union of circles of different sizes. Test problems 2-5 are  annular containers, which have a  single inner circular prohibited area. Test problem 6 has prohibited areas that are disconnected.
Note here that given Table~\ref{Table:coordinates} future
workers will be able to compare their results against ours 
for the same set of test problems as we have used.
For each of these six test problems we considered instances with $n=10,20,...,100$ circles. 
Each value of $n$ considered was solved independently, 
in particular we do not use any information from the 
 solution for packing $n$ circles (as found by our heuristic) 
to assist in creating the solution for 
 packing  $n+10$ circles (or vice-versa).

\begin{table}[!htb] 
\begin{center}
\caption{Centre coordinates and radii of the prohibited areas}
\vspace{.2cm}
\begin{tabular}{ccccc}
Test problem  & Circle & Centre x-coordinate & Centre y-coordinate & Radius \\ \hline
 &&\\
1 & 1 & -9.8/15  &  6.3/15 &  1.3/15\\
 & 2 &  -8.2/15  &  6.5/15 &   1.1/15\\
  & 3 & -6.5/15    &6.5/15   & 1.4/15\\
 & 4 &   -6.0/15   & 5.5/15   & 1.3/15\\
  & 5 &  -5.0/15   & 4.0/15   & 1.25/15\\
  & 6 & -4.5/15   & 2.0/15   & 1.5/15\\
 & 7 &  -4.0/15   & 0.5/15   & 1.25/15\\
 & 8 &  -3.5/15   &-1.0/15   & 1.35/15\\
 &  9 &   -2.8/15   & 2.7/15    &1.2/15\\
 & 10 &   -6.3/15   & 1.7/15    &1.2/15\\
 & 11 &  -1.3/15   & 3.7/15    &1.2/15\\	
	 \\ \hline
\\
 2 & 1 & 0 & 0 & 1/10.5 \\
	 \\ \hline
\\
 3 & 1 & 0 & (1/10.5)-1 & 1/10.5\\ 
	 \\ \hline
\\
 4 & 1 & 0 & 0 & 10.25/17.5 \\
	 \\ \hline

\\
 5 & 1 & 0 & (10.25/17.5)-1 & 10.25/17.5\\ 
 \\ \hline
\\
6 & 1 & 0 	&	(1/10.5)-1 	& 1/10.5 \\
& 2 & 0	&	1-(1/10.5) 	& 1/10.5 \\
& 3 & (1/10.5)-1 &	0	&	1/10.5 \\
& 4 & 1-(1/10.5)	& 0	&	1/10.5 \\
	 \\ \hline
	\end{tabular}
	\label{Table:coordinates}
\end{center}
\end{table}

In Table~\ref{Table:fa3} we present the results obtained for test problems 2-5. These are all cases for which we have just  a single circular prohibited area.
In that table we show, for each value of $n$ considered, the best (maximum) radius found over all replications, together with the total time (in seconds) for all replications.
In general, considering that for each value of $n$ we have to solve 2000 nonlinear programs, we can say that the total time required  appears reasonable.


In~Figures~2(a)-2(d)
we show  visually the results for test problems 2-5 with $n=50$; and in Figures~3(a)-3(d) the results for $n=100$. These appear very reasonable, with the circles being closely packed around the prohibited area.

As can be seen from Table~\ref{Table:coordinates} test problems 2 and 3 both have a single prohibited area, namely one
 circle with radius (1/10.5).  The difference between them relates to where in the unit circle container this single prohibited area is positioned 
(cf~Figures 2(a) and 2(b); Figures 3(a) and 3(b)). Since the area of the unit circle  container that can be used 
(i.e.~the non-prohibited area) is hence the same in both cases it might seem reasonable to suppose that, when we attempt to maximise the radius of $n$ identical packed circles, the best radius found and the time required would be similar for this pair of test problems
(for each value of $n$).

For test problems 2 and 3 the average absolute difference between the best radii  is 0.00216184, so very small. As a percentage of the average radius for test problem 2 this difference is only  1.58\%, and as a percentage of the average radius for test problem 3 this difference is only  1.56\%. However the average time required for test problem 2 is 127.07 seconds, whereas for test problem 3 it is much larger, 303.31 seconds. These results therefore indicate that (for our FSS algorithm) the position of the prohibited area has a significant effect on computation time.

In a similar fashion,  as can be seen from 
Table~\ref{Table:coordinates}, test problems 4 and 5 both have a single prohibited area, namely one
 circle with radius (10.25/17.5).  The difference between them relates to where in the unit circle container this single prohibited area is positioned 
(cf~Figures 2(c) and 2(d); Figures 3(c) and 3(d)). So again it might seem reasonable to suppose that, when we attempt to maximise the radius of $n$ identical packed circles, the best radius found and the time required would be similar for this pair of test problems (for each value of $n$). 

For these test problems the average absolute difference between the best radii  is 0.00254764, so again very small. As a percentage of the average radius for test problem 4 this difference is only  
2.37\%, and as a percentage of the average radius for test problem 5 this difference is only  2.35\%. However the average time required for test problem 4 is 88.72 seconds, whereas for test problem 5 it is much larger, 499.56 seconds.  These results therefore again indicate that (for our FSS algorithm) the position of the prohibited area has a significant effect on computation time.

To investigate how our FSS approach performed when we had more than one circular prohibited area we used test problem 1, but varied the number of circles $|F|$ in the prohibited area. 
The results can be seen in Table~\ref{Table:fa5} and Table~\ref{Table:fa6}. Here we again show the best (maximum) radius found and the total time taken (in seconds). The results in Table~\ref{Table:fa5} for  $|F|=6$, for example, correspond to using just the first six circles associated with that test problem (the size and position of these  circles being given in Table~\ref{Table:coordinates}).

In~Figures~2(e)-2(l) 
we show visually the results for $n=50$ and test problem 1 as $|F|$ varies. Figures~3(e)-3(l)  show the results for $n=100$. Note from these figures that the prohibited area for test problem
 1,  composed of $|F|$ circles, is overall irregular in shape (but always connected).

Table~\ref{Table:fa5} and Table~\ref{Table:fa6} indicate that the average computation time increases as $|F|$ increases, as we might expect. However the increase is not especially marked, the minimum average computation time in those tables is 392.21 seconds associated with $|F|=5$, whilst the maximum average computation time is 651.76 seconds associated with $|F|=10$.

Visually it is clear from comparing the results for test problems 2-5 with those for test problem 1 that the results for test problem 1 exhibit more unused space, principally around the irregular prohibited area. However the portion of the circular container that does not involve the prohibited area is relatively closely packed with circles.

The reason, we believe, why these results exhibit unused space around the  irregular prohibited area is that as $n$ increases the predominant factor in maximising radius becomes the packing of the circles in the region away from the prohibited area. For example consider 
Figure~3(l) with $n=100$ circles. First recall that the number of circles is fixed, so we cannot fill the unused space seen in that figure with extra circles of the same radius as those already shown in Figure~3(l). Increasing, by even a small amount, the (equal) radii of the circles shown, will clearly be difficult for the closely packed circles away from the prohibited region
(even allowing for the fact that all circles can be resized and repositioned).

As mentioned above, our computational experience has been that introducing 
Equations~(\ref{e14})-(\ref{e15}) does significantly improve the quality of the results that we obtain. As an illustration of this 
we solved all problems (from $n=10$ to $n=30$, test problems 1-5)  both with, and without,  Equations~(\ref{e14})-(\ref{e15}). Overall, when solving with Equations~(\ref{e14})-(\ref{e15}) included, the average computation time was 33.5\% \textbf{\emph{lower}} and the average best radius 12.8\% \textbf{\emph{higher}}, as compared with the situation when these equations were not included.  
In other words we (on average) obtain much better results, both with respect to computation time and with respect to the value of the best radius found, when we include Equations~(\ref{e14})-(\ref{e15}).
This vindicates, we believe, our algorithmic design decision to include Equations~(\ref{e14})-(\ref{e15}).

With regard to test problem 6 Table~\ref{Table:dis} shows the results obtained.
In~Figures~4(a)-4(b)
we show  visually the results for that test problem with $n=50$ and $n=100$. It can be seen that in this particular instance we have four disconnected prohibited areas. Again the results appear  reasonable, with the circles being closely packed around the prohibited areas.

\section{Conclusions} \label{Sec:Conclusions}

In this paper we have considered the problem of packing a fixed number of identical circles inside the unit circle container, where the packing was complicated by the presence of fixed size circular prohibited areas. The objective we adopted was  to 
maximise the radius of the identical circles.

We considered six different test problems, one with a single prohibited area being the union of circles of different sizes, four being 
different annular containers and one where the  prohibited areas were disconnected.

We presented  a heuristic for the problem based upon formulation space search, a new and emerging metaheuristic. Computational results were given for problems involving the packing of up to 100 circles.

\clearpage
\newpage

\begin{table}[!htb] 
{\scriptsize
\begin{center}
\caption{FSS results  for test problems 2-5}
\vspace{.2cm}
\begin{tabular}{ccccccccc} \hline
	&		\multicolumn{2}{c}{Test problem 2}	&	\multicolumn{2}{c}{Test problem 3}	&	\multicolumn{2}{c}{Test problem 4}	&	\multicolumn{2}{c}{Test problem 5}	\\	
n	&	Best radius	&	Total time (s) &	Best radius	&	Total time (s)	&	Best radius	&	Total time (s)	
&	Best radius	&	Total time (s)	\\
\hline	
10		&	0.25060817	&	18.25	&	0.26225892	&	64.95
&	0.20714286	&	12.73	&	0.20620478	&	193.20
	\\
20		&	0.19039215	&	28.62	&	0.19522401	&	72.58	
&	0.14044117	&	20.22	&	0.14956309	&	229.56
\\
30		&	0.15919784	&	46.30	&	0.15979391	&	136.80	
&	0.11926162	&	29.70	&	0.12434830	&	340.39
\\
40		&	0.13742599	&	60.15	&	0.13930887	&	170.70	
&	0.11078813	&	44.56	&	0.10879452	&	246.51
\\
50		&	0.12471293	&	99.35	&	0.12517615	&	210.50	
&	0.10055446	&	66.65	&	0.09792755	&	377.35
\\
60		&	0.11545614	&	103.10	&	0.11545599	&	261.50	
&	0.08859160	&	81.85	&	0.08984271	&	304.23
\\
70		&	0.10533517	&	156.50	&	0.10605514	&	309.20	
&	0.08215214	&	104.03	&	0.08377036	&	733.21
\\
80		&	0.09972555	&	205.90	&	0.09916911	&	454.00	
&	0.07827693	&	139.87	&	0.07865683	&	663.07
\\
90		&	0.09460328	&	237.40	&	0.09370228	&	587.50	
&	0.07473140	&	170.12	&	0.07429912	&	922.83
\\
100		&	0.08899120	&	315.10	&	0.08900724	&	765.40	
&	0.07239415	&	217.44	&	0.07036645	&	985.28
\\
\hline
\end{tabular}
\label{Table:fa3}
\end{center}
}
\end{table}

\vspace{.2cm}
\begin{table}[!htb] 
{\scriptsize
\begin{center}
\caption{FSS results  for test problem 1 with $|F|$ varying from 4 to 7}
\vspace{.2cm}
\begin{tabular}{ccccccccc} \hline
	&	\multicolumn{2}{c}{$|F| = 4$}
	&	\multicolumn{2}{c}{$|F| = 5$} &	\multicolumn{2}{c}{$|F| = 6$}&	\multicolumn{2}{c}{$|F| = 7$} 	\\	
	n	&	Best radius	&	Total time (s)&	Best radius	&	Total time (s)&	Best radius	&	Total time (s)
&	Best radius	&	Total time (s)	\\
\hline	
10 &	0.25725385	&	262.81	&	0.25725385	&	232.91	&	0.25725385	&	271.90	&	0.25725385	&	284.46	\\
20	&	0.18800266	&	261.32 &	0.18590633	&	289.80	&	0.18531662	&	247.60	&	0.18531515	&	446.41	\\
30		&	0.15561785	&	292.92 &	0.15589953	&	296.77	&	0.15398348	&	309.10	&	0.15252655	&	489.18	\\
40		&	0.13597134	&	250.93 &	0.13494222	&	210.46	&	0.13301332	&	237.80	&	0.13273445	&	404.74	\\
50	&	0.12225672	&	254.90 &	0.12083527	&	267.85	&	0.12021850	&	285.80	&	0.12068548	&	357.89	\\
60	&	0.11019235	&	316.74 &	0.11030502	&	296.81	&	0.10935882	&	433.90	&	0.10898846	&	429.68	\\
70	&	0.10343522	&	420.70 &	0.10268332	&	450.44	&	0.10153295	&	461.80	&	0.10066682	&	500.23	\\
80	&	0.09611034	&	543.14 &	0.09631517	&	605.45	&	0.09448850	&	625.00	&	0.09262220	&	703.18	\\
90	&	0.09201990	&	532.95 &	0.09149489	&	544.22	&	0.08825820	&	626.40	&	0.08757232	&	787.02	\\
100	&	0.08683729	&	844.48 &	0.08651290	&	727.39	&	0.08464148	&	983.10	&	0.08446274	&	980.21	\\
  \hline
\end{tabular}
\label{Table:fa5}
\end{center}
}
\end{table}

\vspace{.2cm}
\begin{table}[!htb] 
{\scriptsize
\begin{center}
\caption{FSS results   for test problem 1 with $|F|$ varying from 8 to 11}
\vspace{.2cm}
\begin{tabular}{ccccccccc} \hline
	&	\multicolumn{2}{c}{$|F| = 8$}&	\multicolumn{2}{c}{$|F| = 9$}&	\multicolumn{2}{c}{$|F| = 10$} 			&	\multicolumn{2}{c}{$|F| = 11$} \\	
n	&	Best radius	&	Total time (s)&	Best radius	&	Total time (s)&	Best radius	&	Total time (s)	
&	Best radius	&	Total time (s)	\\
 \hline	
10	&	0.25582765	&	277.31	&	0.25297675	&	245.78	&	0.25213340	&	280.58	
&	0.24958389	&	249.10
\\
20	&	0.18395004	&	318.46	&	0.18278297	&	379.61	&	0.18145427	&	402.61	
&	0.17857572	&	390.90
\\
30	&	0.15253463	&	473.42	&	0.15100605	&	502.87	&	0.15099462	&	632.30	
&	0.14944660	&	685.50

\\
40	&	0.13238624	&	322.10	&	0.13150801	&	316.52	&	0.12939828	&	415.40	
&	0.12945512	&	379.90

\\
50	&	0.11841380	&	363.05	&	0.11776737	&	372.58	&	0.11879430	&	557.99	
&	0.11775790	&	466.00

\\
60	&	0.10685984	&	343.93	&	0.10841268	&	412.69	&	0.10596339	&	518.94	
&	0.10633808	&	482.20

\\
70	&	0.09940636	&	430.45	&	0.09976116	&	511.19	&	0.09981981	&	688.71	
&	0.09842435	&	635.50
\\
80	&	0.09187960	&	520.99	&	0.09130728	&	562.96	&	0.09215045	&	811.07	
&	0.09195809	&	774.40

\\
90	&	0.08838363	&	670.41	&	0.08841618	&	964.66	&	0.08515947	&	1145.45	
&	0.08674741	&	785.40

\\
100	&	0.08353399	&	770.81	&	0.08299581	&	1169.21	&	0.08350246	&	1064.56	
&	0.08286614	&	1184.00

\\
 \hline
\end{tabular}
\label{Table:fa6}
\end{center}
}
\end{table}

\vspace{.2cm}
\begin{table}[!htb] 
{\scriptsize
\begin{center}
\caption{FSS results  for test problem 6}
\vspace{.2cm}
\begin{tabular}{ccc} \hline
n	&	Best radius	&	Total time (s) \\
\hline	
 
10	& 	0.26018588	&	122.90	\\	   
20	& 	0.18808326	& 	133.40	\\	   
30	& 	0.15423705	& 	232.40	\\	   
40	& 	0.13318076	& 	194.00	\\	   
50	& 	0.11909887	& 	299.60	\\	   
60	& 	0.10972674	& 	440.00	\\	   
70	& 	0.09834932	& 	607.90	\\	   
80	& 	0.09429742	& 	928.70	\\	   
90	& 	0.08940553	& 	737.40	\\	   
100	& 	0.08427842	& 	962.60	\\	 

\hline
\end{tabular}
\label{Table:dis}
\end{center}
}
\end{table}

\begin{landscape}
\begin{figure}[!htbp] 
\centering   
\caption{FSS results for packing 50 identical circles in a circular container for test problems 1-5}
\vspace{-.5cm}
\subfigure[Test problem 2]{\label{fig:50fa2}\includegraphics[scale=0.4]{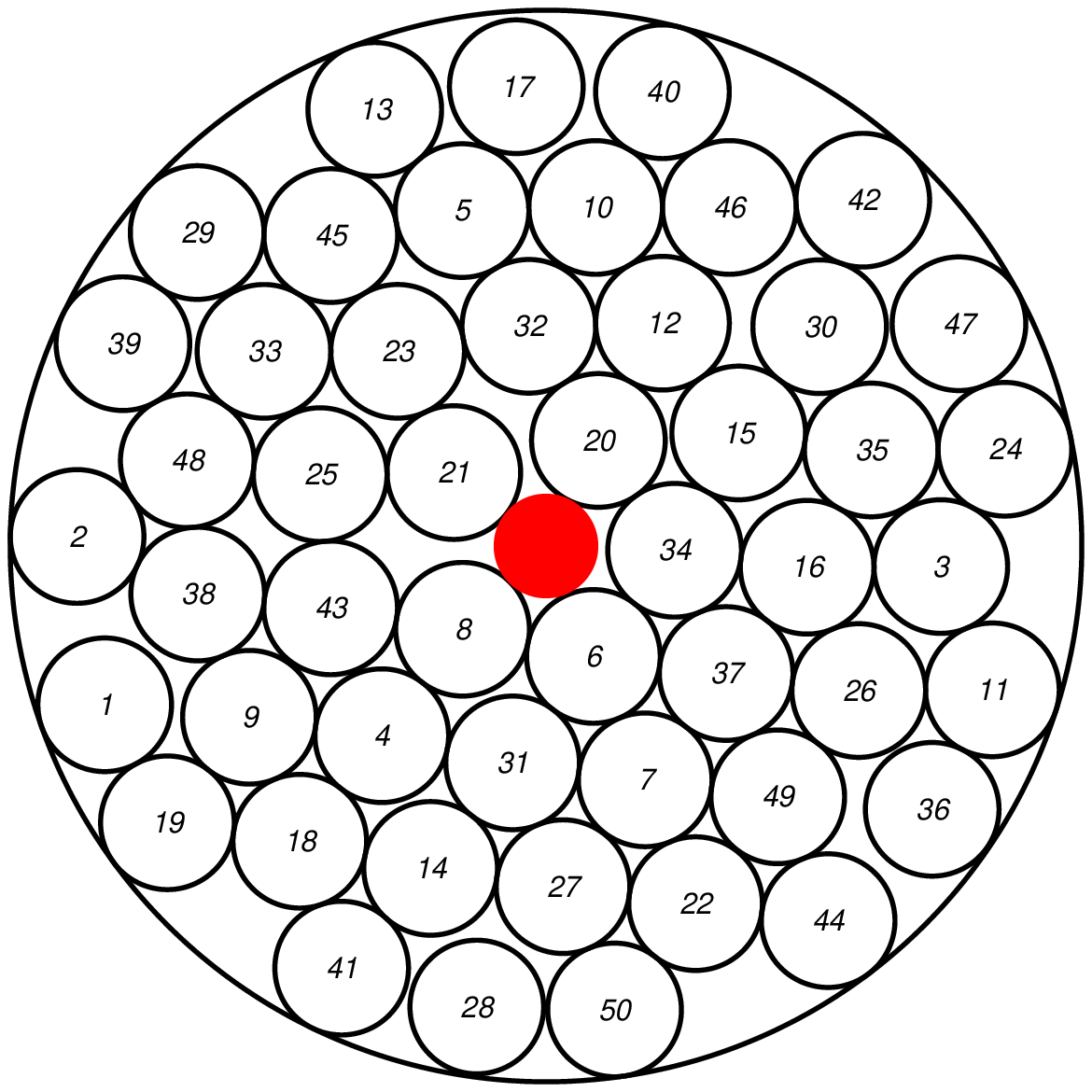}} \hspace{.2cm}
\subfigure[Test problem 3]{\label{fig:50fa3}\includegraphics[scale=0.4]{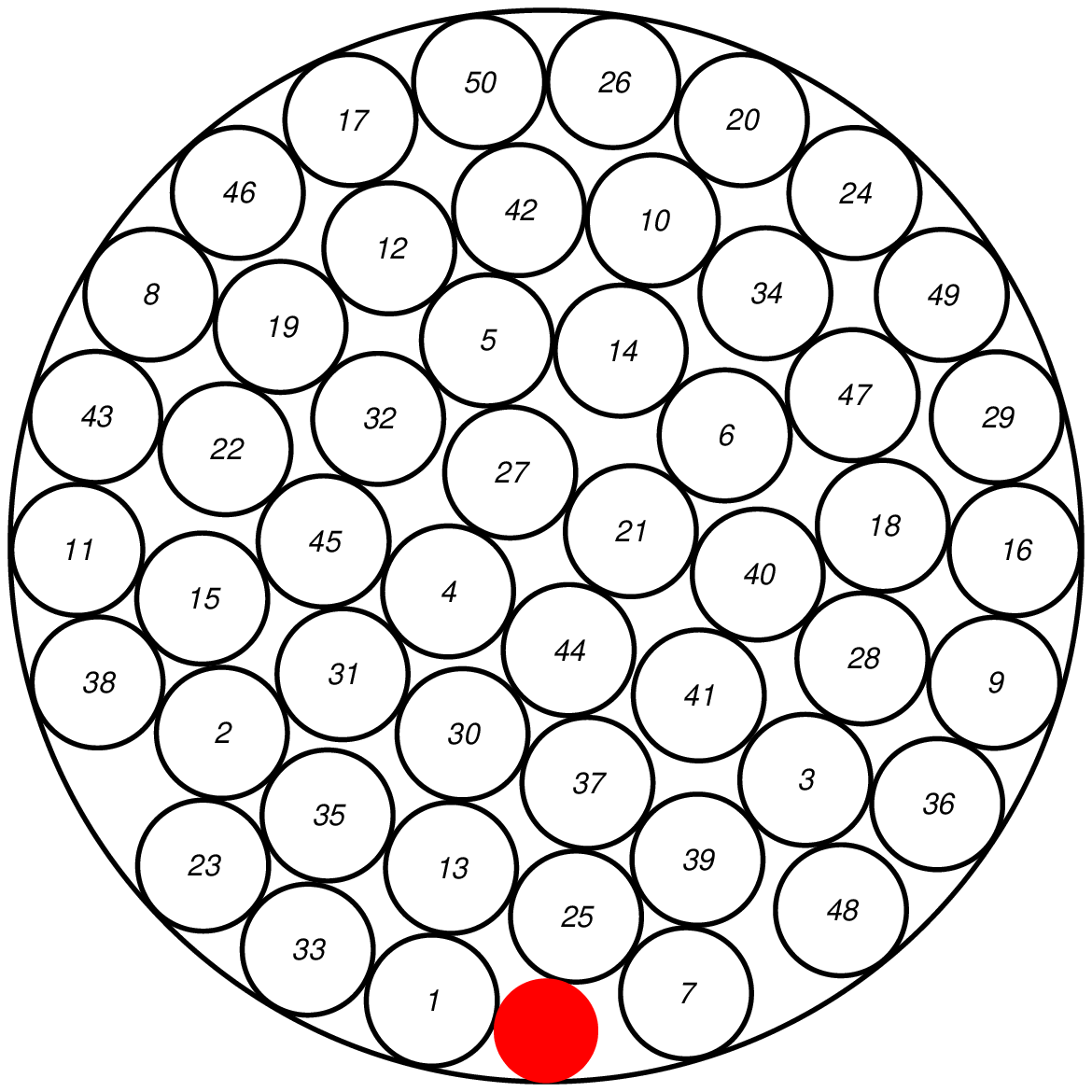}} \hspace{.2cm}
\subfigure[Test problem 4]{\label{fig:50fa4}\includegraphics[scale=0.4]{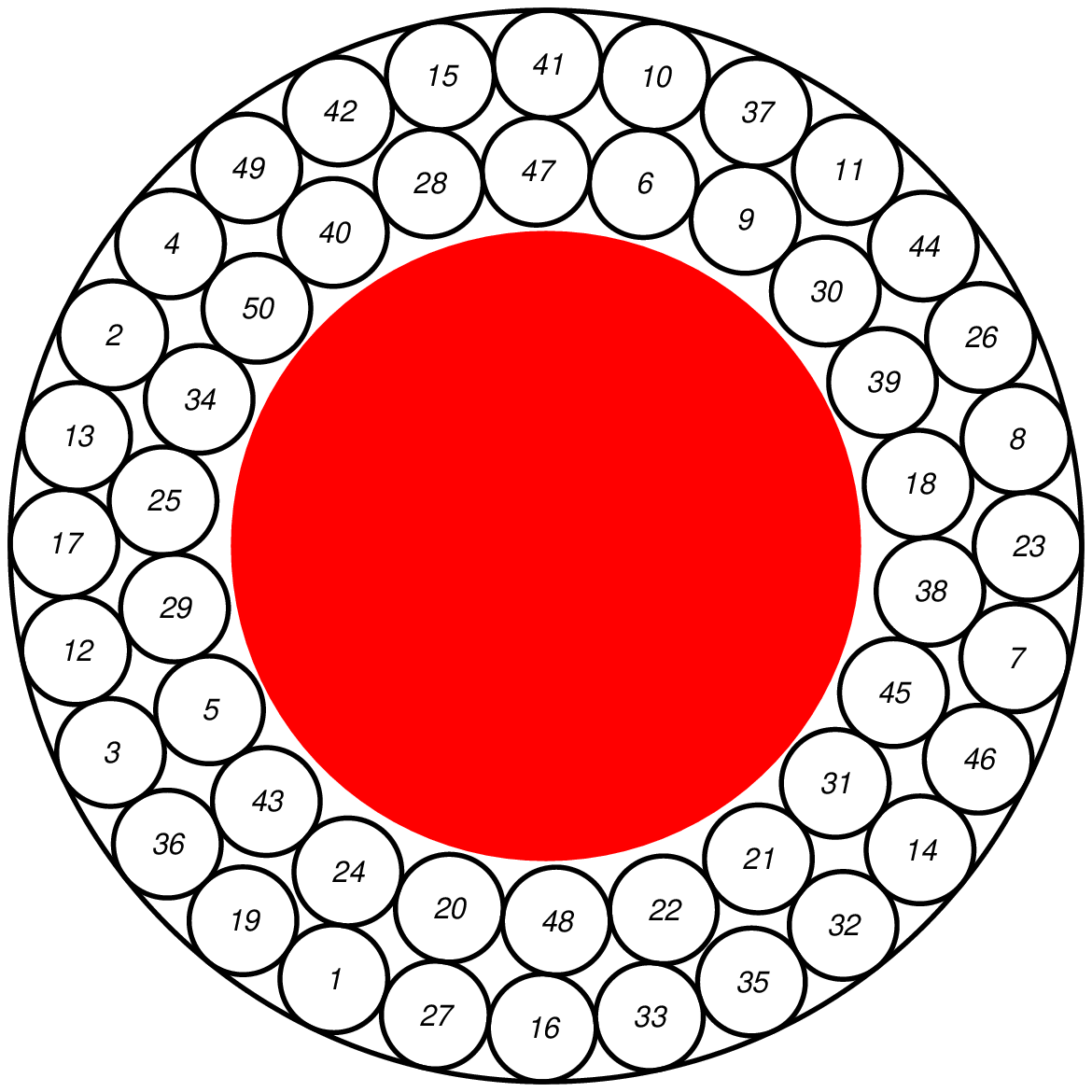}} \hspace{.2cm}
\subfigure[Test problem 5]{\label{fig:50fa5}\includegraphics[scale=0.4]{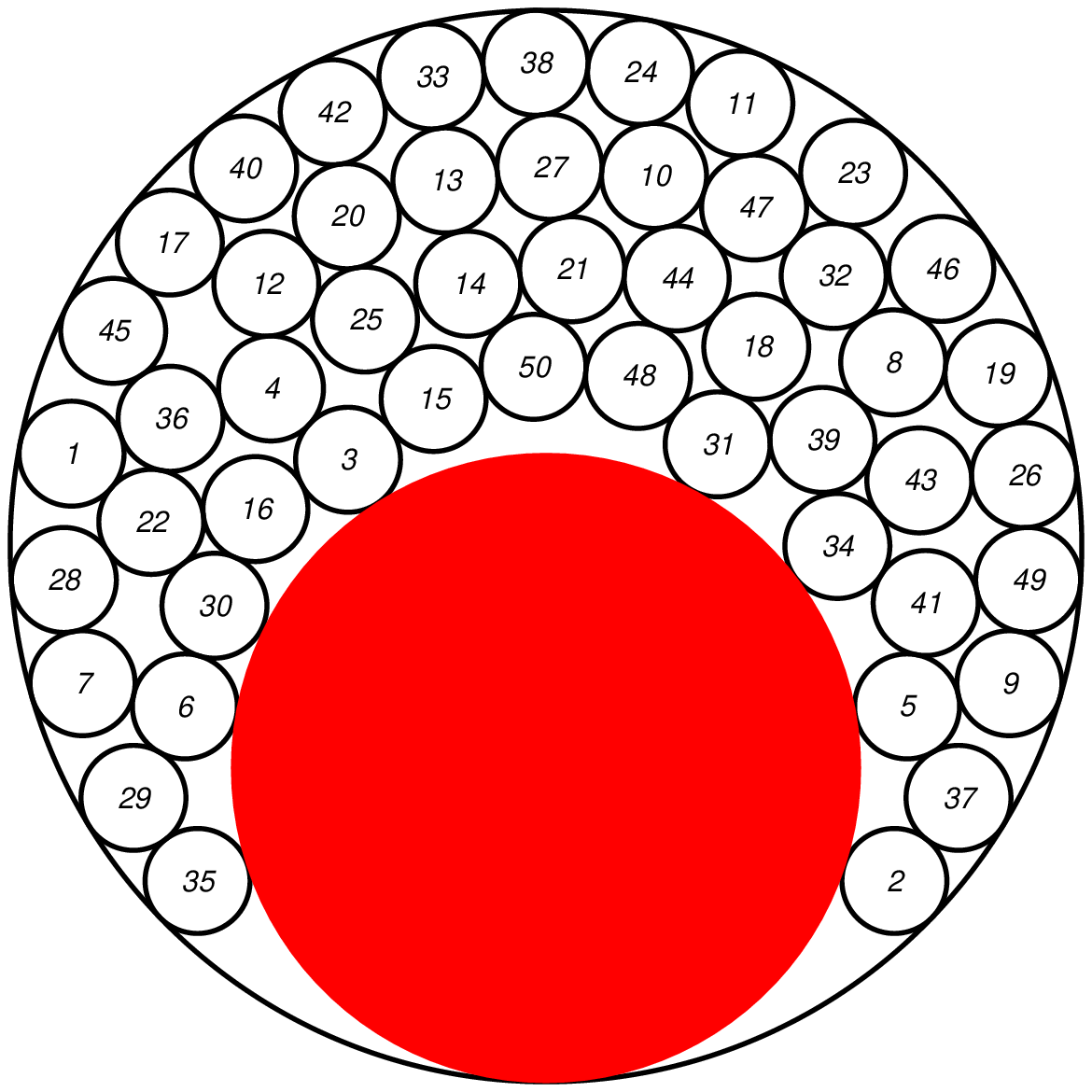}} 
\\
\subfigure[Test problem 1,  $|F|=4$]{\label{fig:50fac4}\includegraphics[scale=0.4]{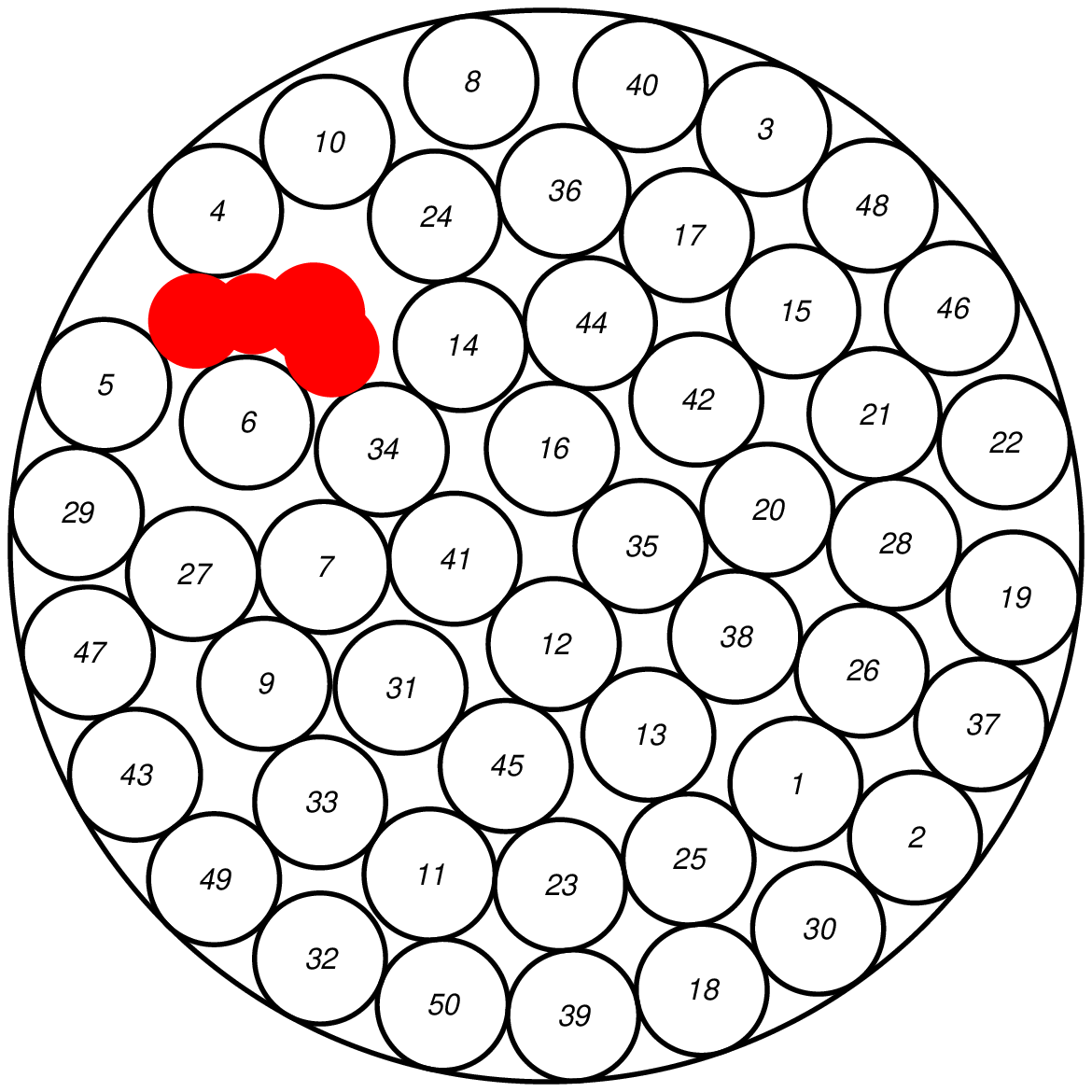}} \hspace{.2cm}
\subfigure[Test problem 1,  $|F|=5$]{\label{fig:50fac5}\includegraphics[scale=0.4]{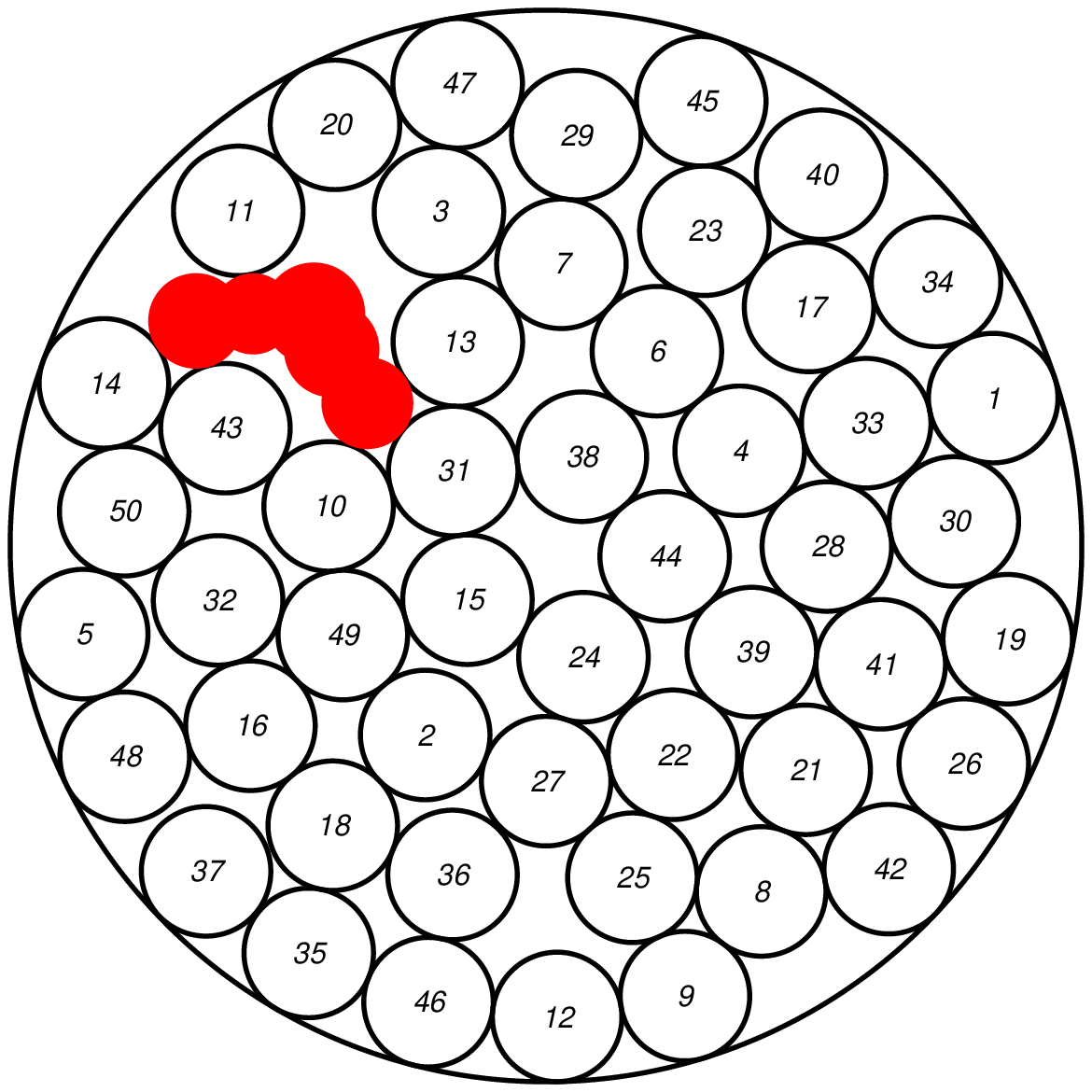}} \hspace{.2cm}
\subfigure[Test problem 1,  $|F|=6$]{\label{fig:50fac6}\includegraphics[scale=0.4]{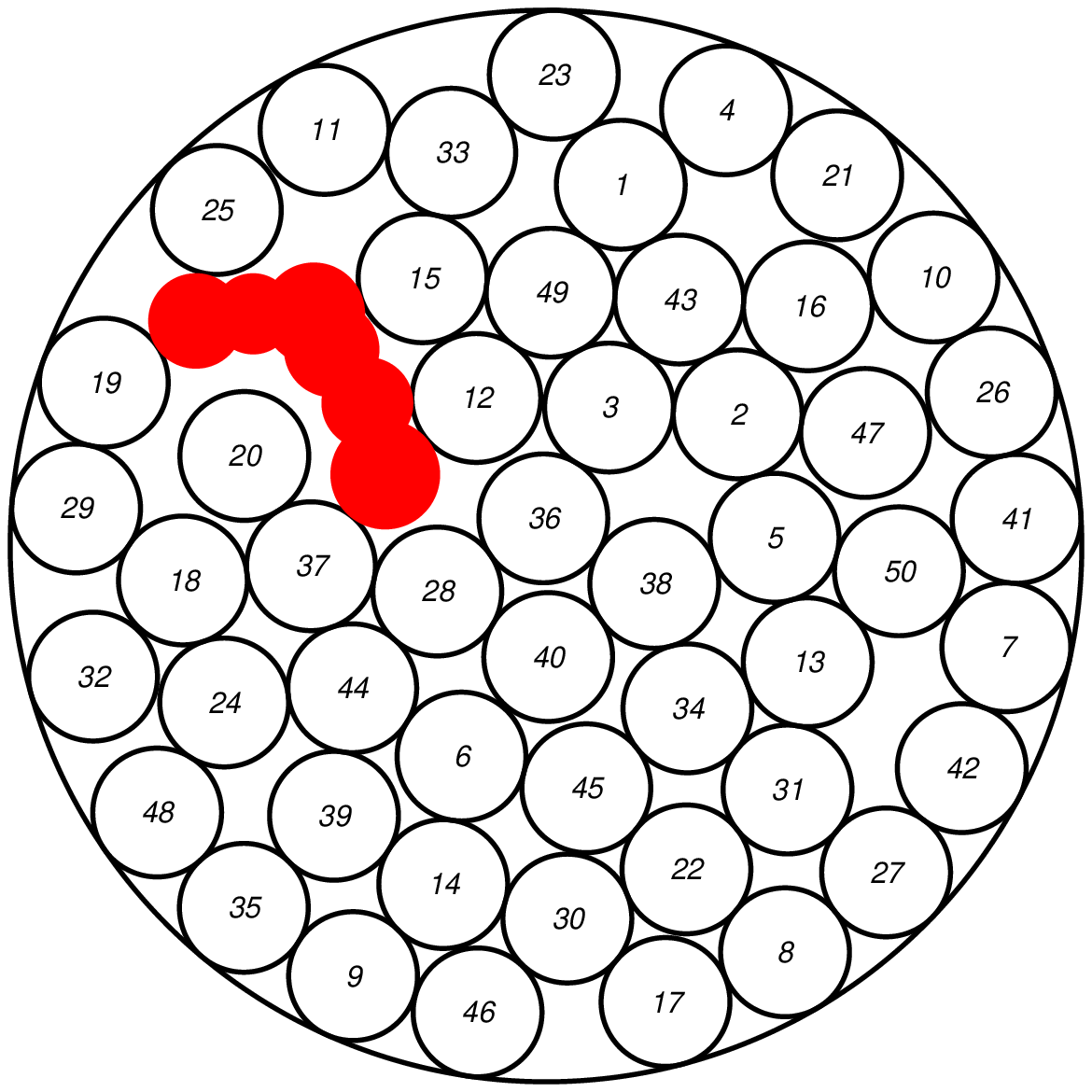}} \hspace{.2cm}
\subfigure[Test problem 1,  $|F|=7$]{\label{fig:50fac7}\includegraphics[scale=0.4]{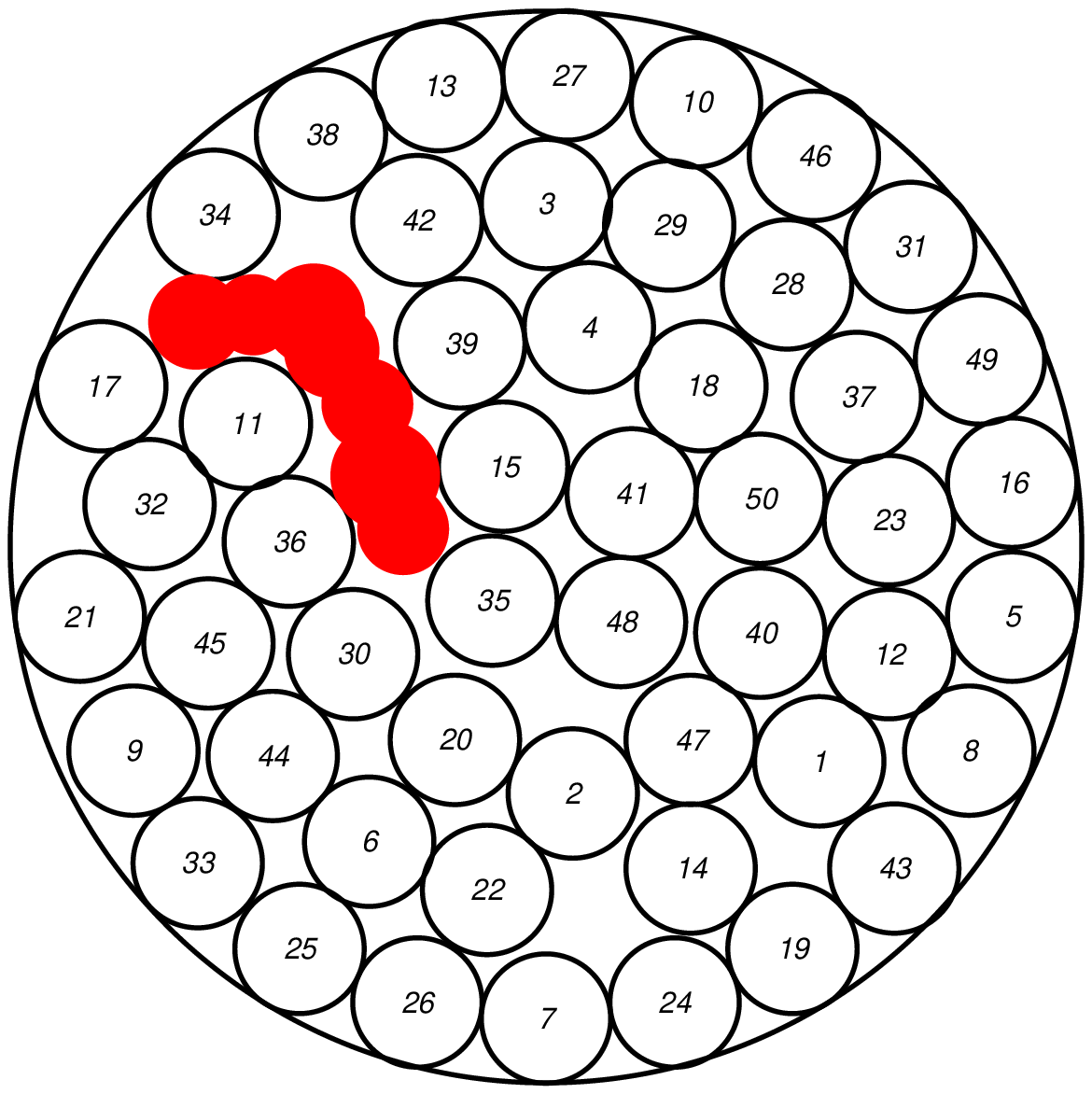}} 
\\
\subfigure[Test problem 1,  $|F|=8$]{\label{fig:50fac8}\includegraphics[scale=0.4]{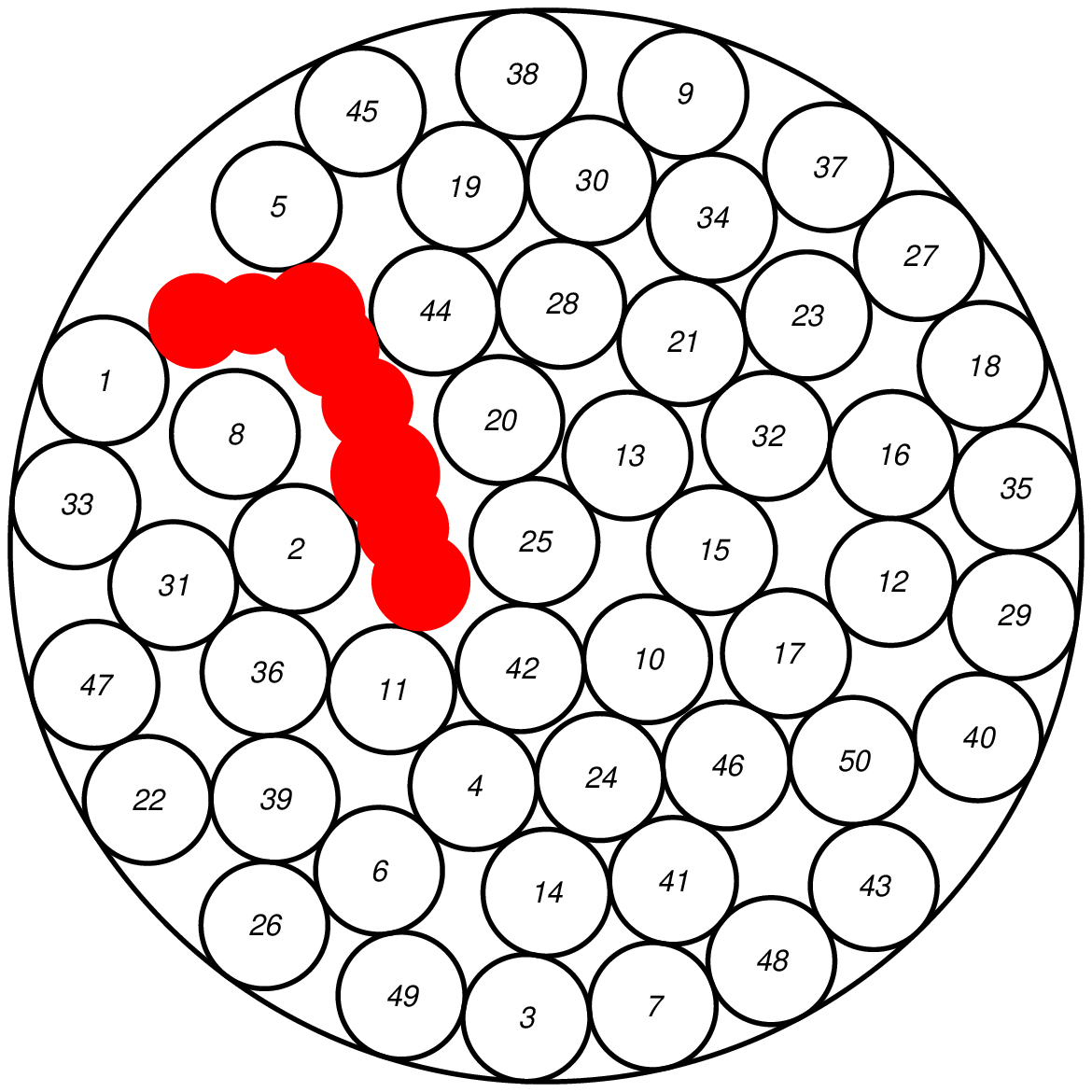}} \hspace{.2cm}
\subfigure[Test problem 1,  $|F|=9$]{\label{fig:50fac9}\includegraphics[scale=0.4]{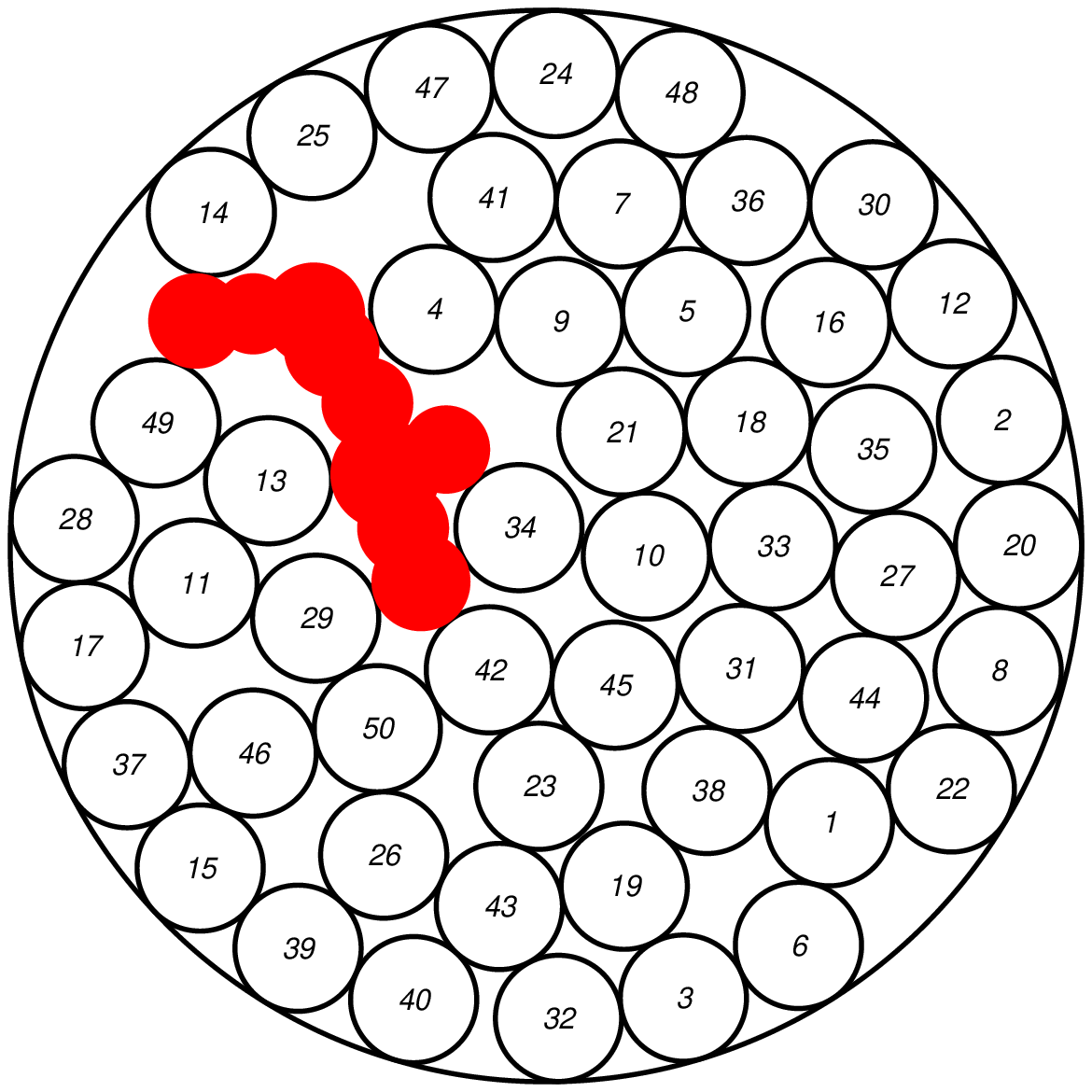}} \hspace{.2cm}
\subfigure[Test problem 1,  $|F|=10$]{\label{fig:50fac10}\includegraphics[scale=0.4]{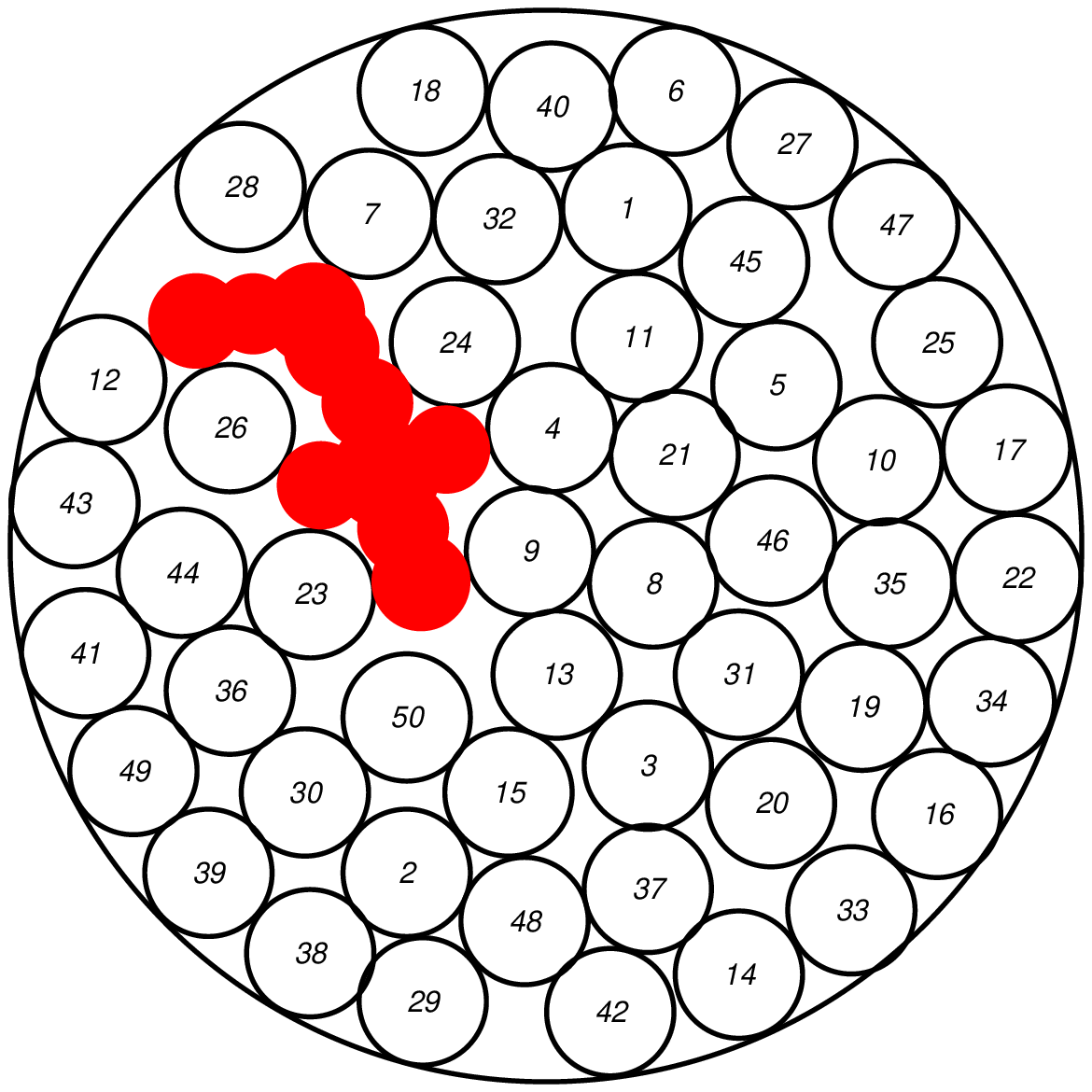}}  \hspace{.2cm}
\subfigure[Test problem 1,  $|F|=11$]{\label{fig:50fa1}\includegraphics[scale=0.4]{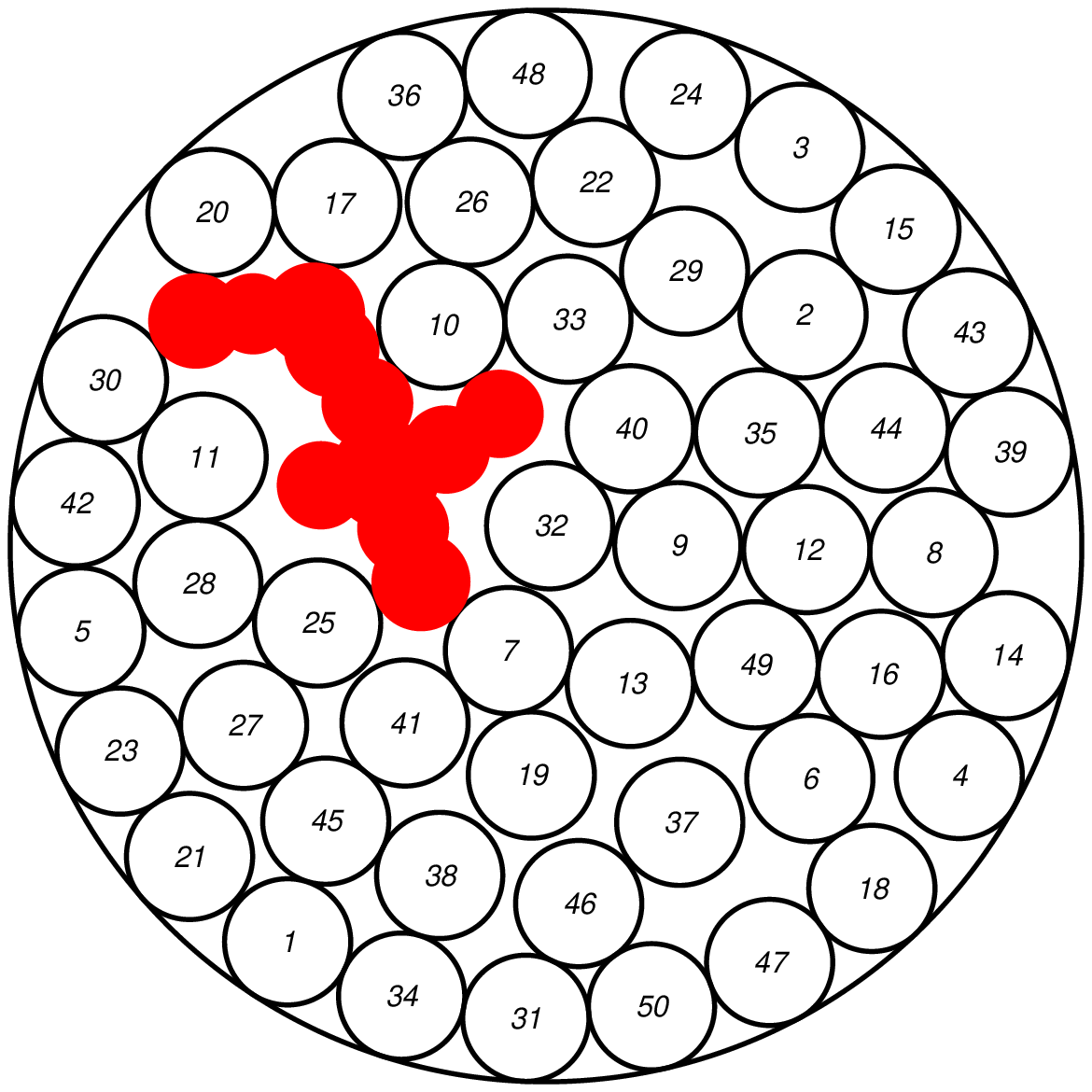}}
\label{Fig:4.50} 
\end{figure}
\end{landscape}

\begin{landscape}
\begin{figure}[!htbp] 
\centering   
\caption{FSS results for packing 100 identical circles in a circular container for test problems 1-5}
\vspace{-.5cm}
 \subfigure[Test problem 2]{\label{fig:100fa2}\includegraphics[scale=0.4]{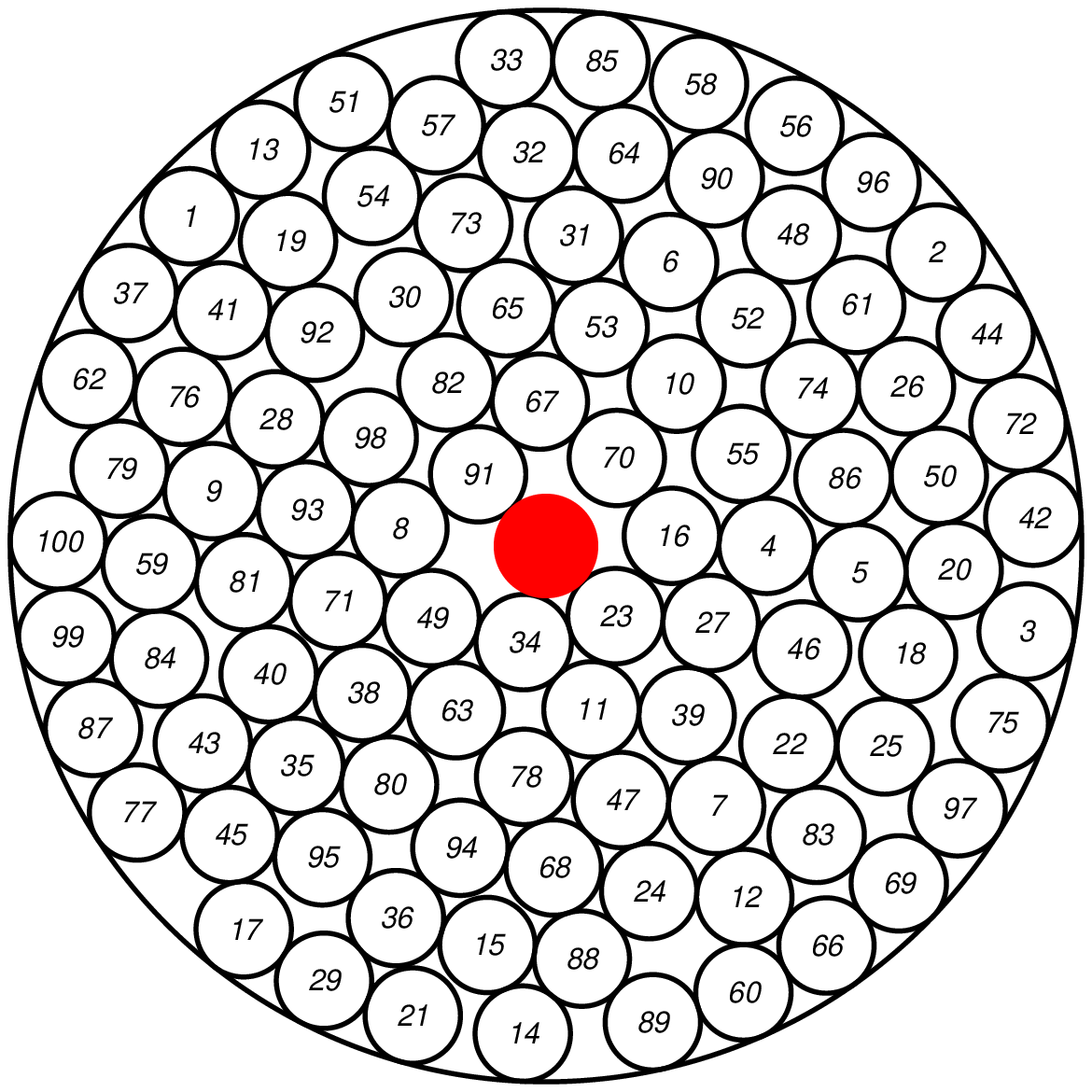}} \hspace{.2cm}
\subfigure[Test problem 3]{\label{fig:100fa3}\includegraphics[scale=0.4]{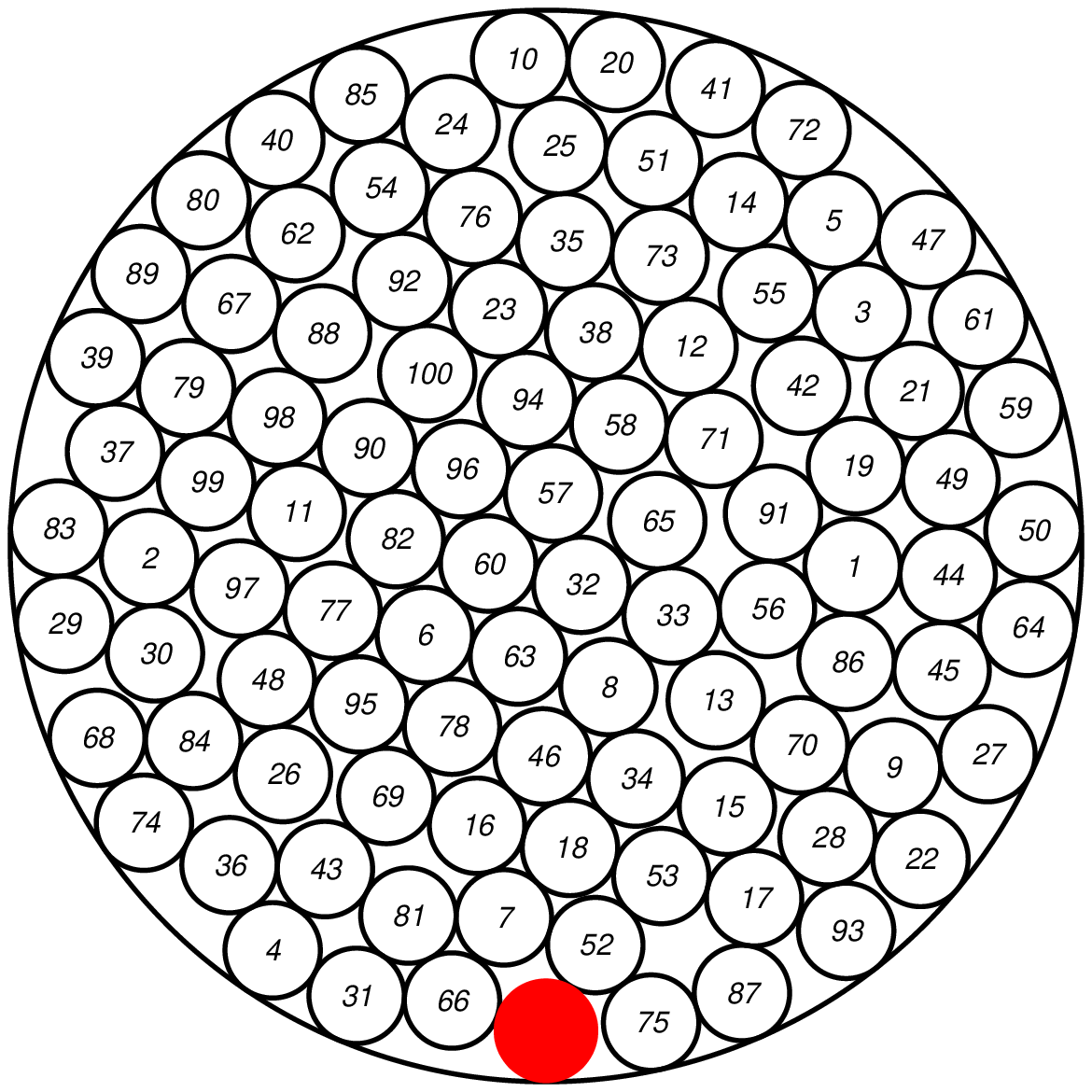}} \hspace{.2cm}
\subfigure[Test problem 4]{\label{fig:100fa4}\includegraphics[scale=0.4]{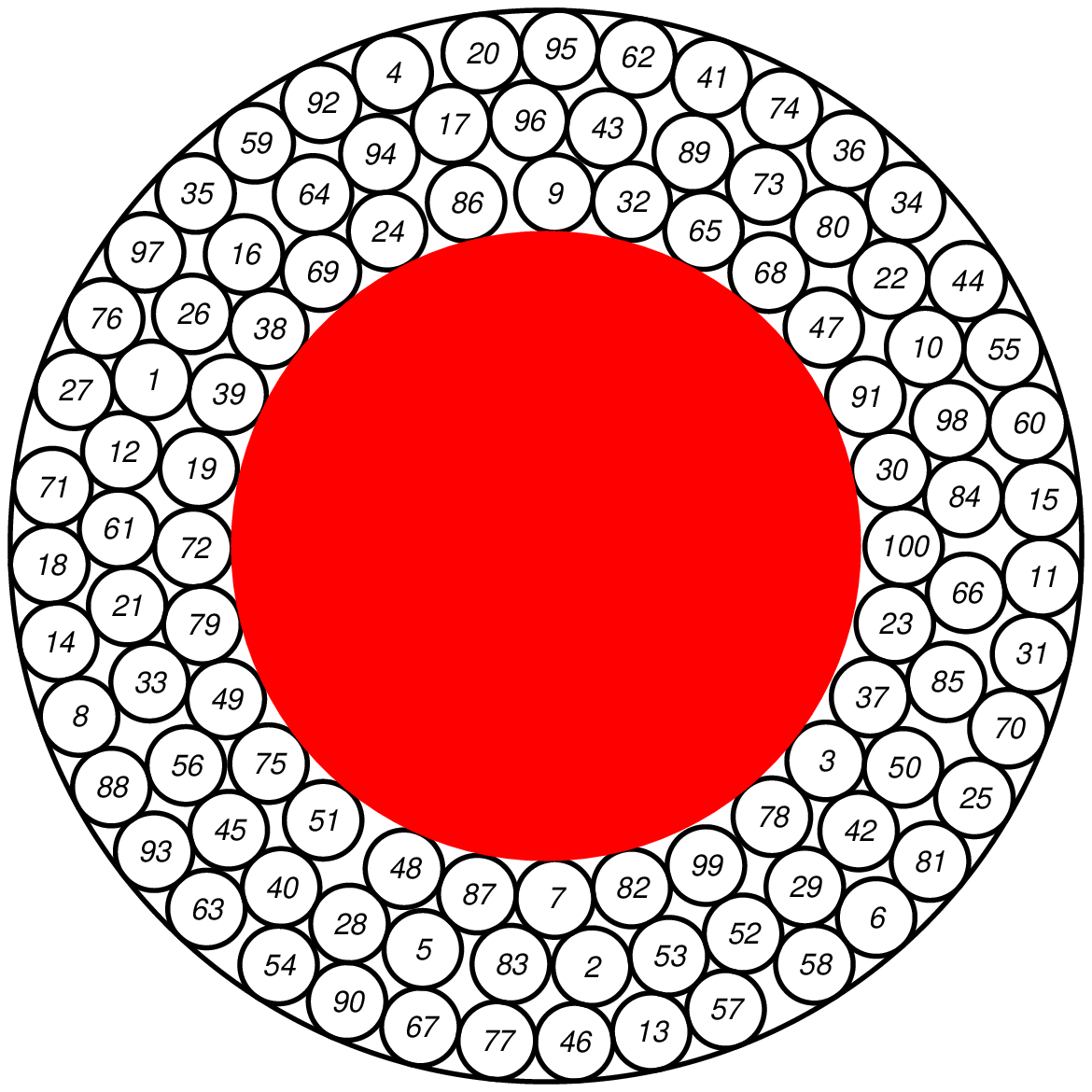}} \hspace{.2cm}
\subfigure[Test problem 5]{\label{fig:100fa5}\includegraphics[scale=0.4]{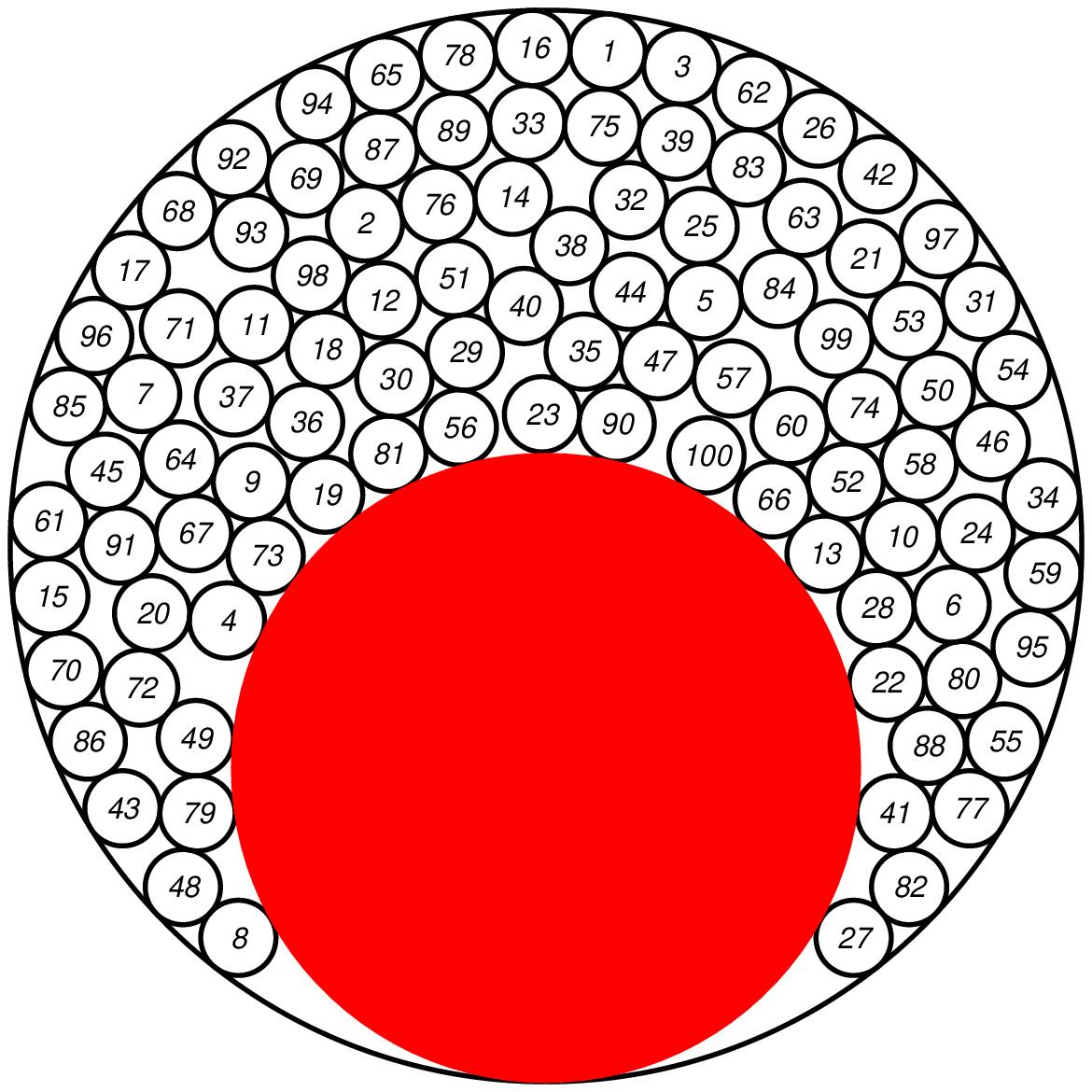}} 
\\
\subfigure[Test problem 1,  $|F|=4$]{\label{fig:100fac4}\includegraphics[scale=0.4]{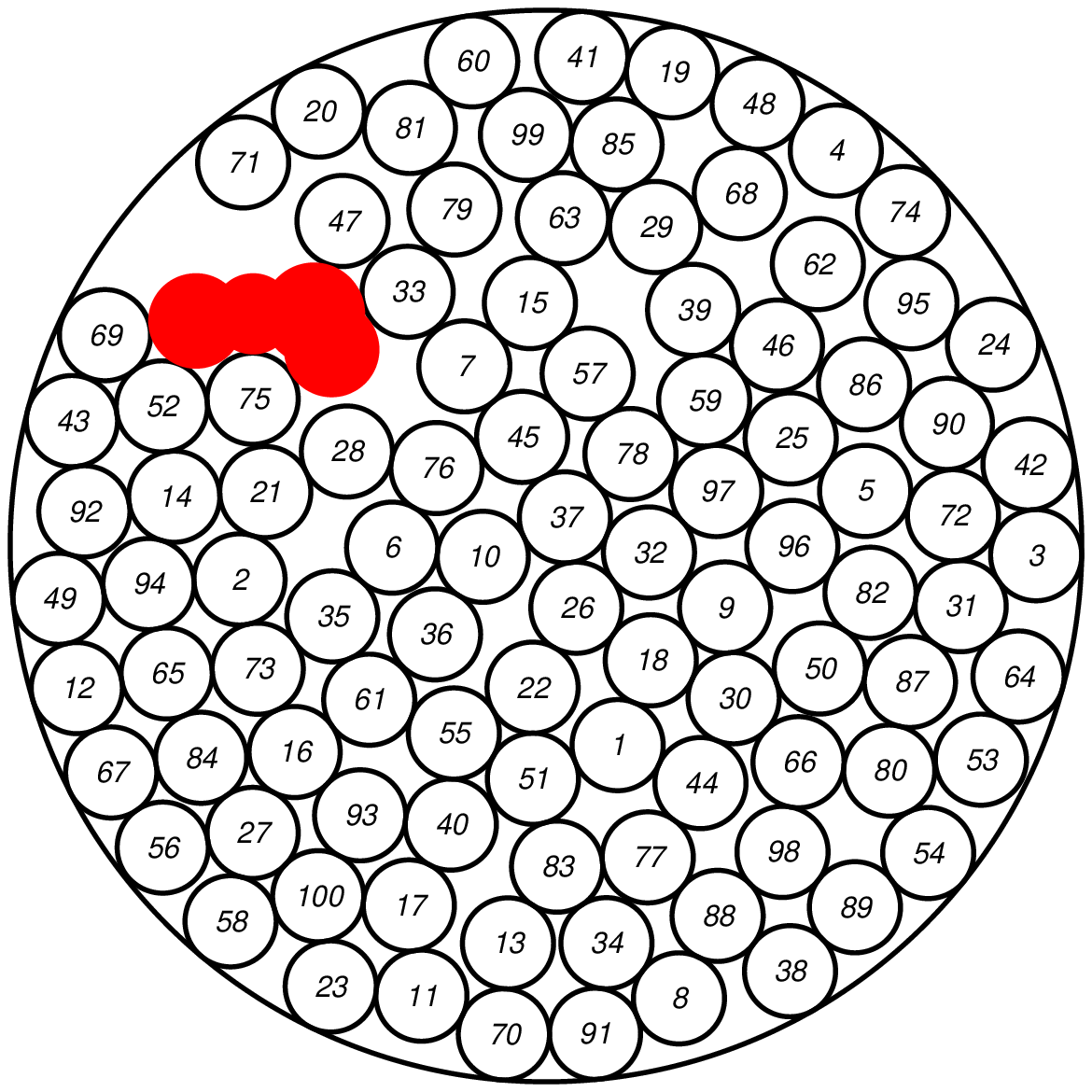}} \hspace{.2cm}
\subfigure[Test problem 1,  $|F|=5$]{\label{fig:100fac5}\includegraphics[scale=0.4]{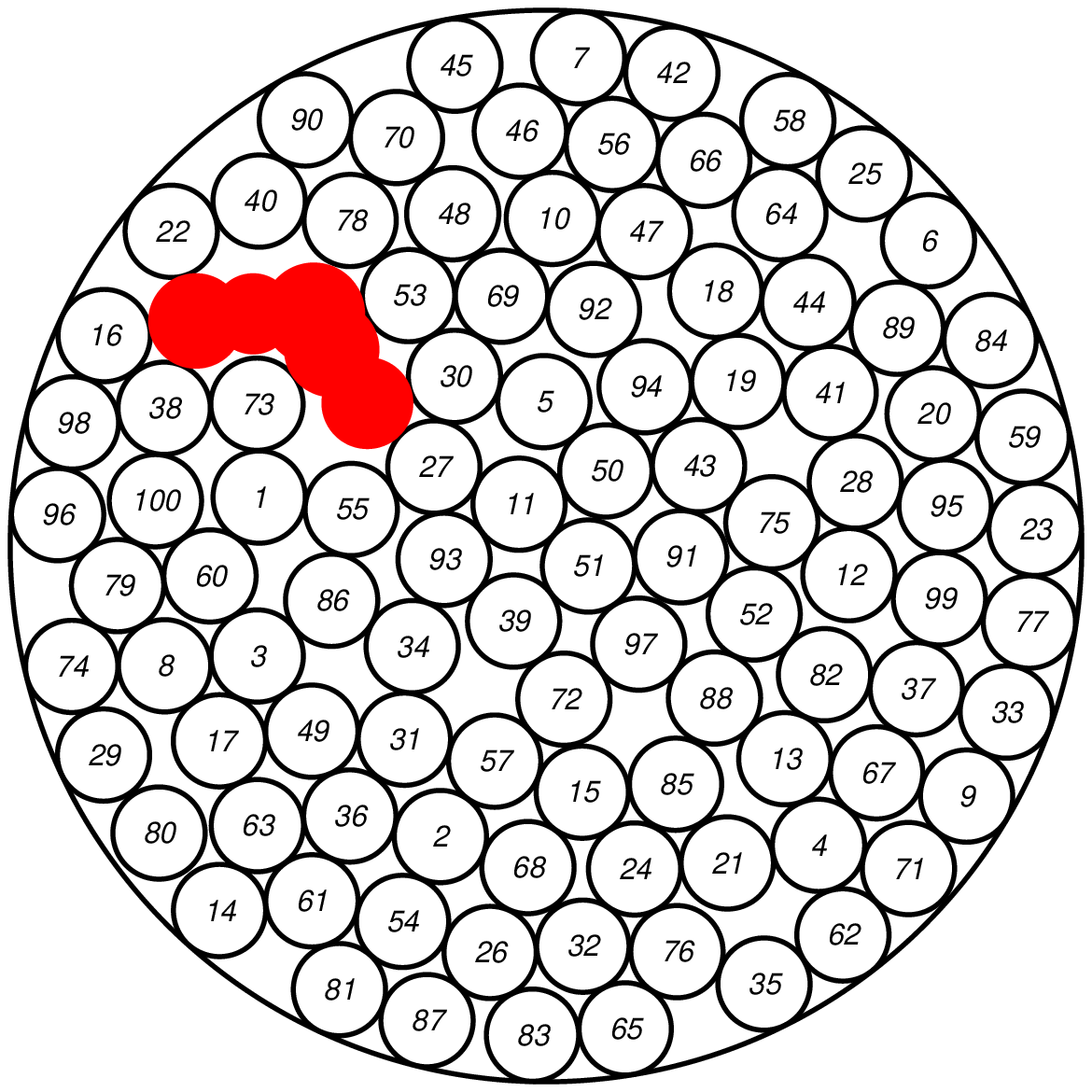}} \hspace{.2cm}
\subfigure[Test problem 1,  $|F|=6$]{\label{fig:100fac6}\includegraphics[scale=0.4]{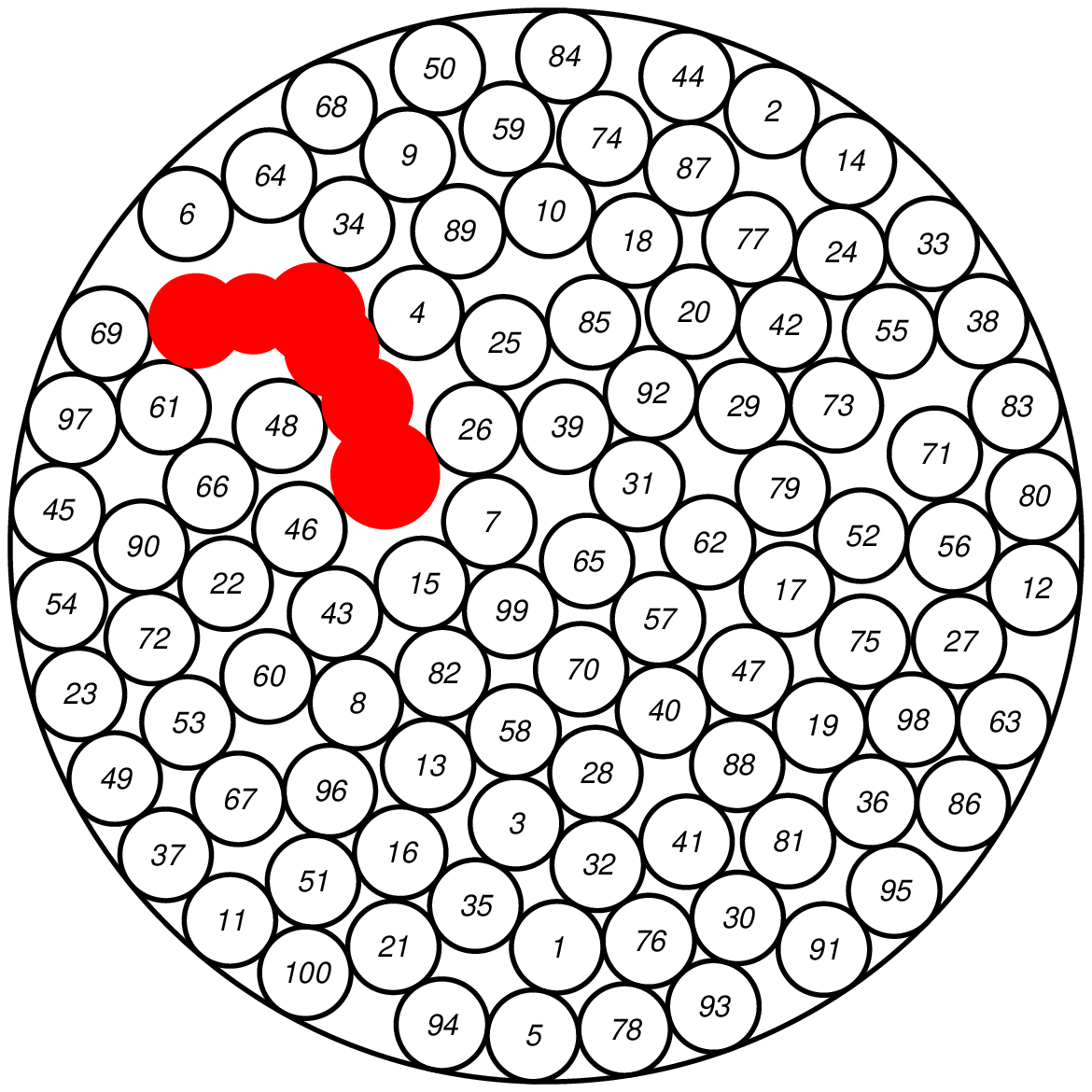}} \hspace{.2cm}
\subfigure[Test problem 1,  $|F|=7$]{\label{fig:100fac7}\includegraphics[scale=0.4]{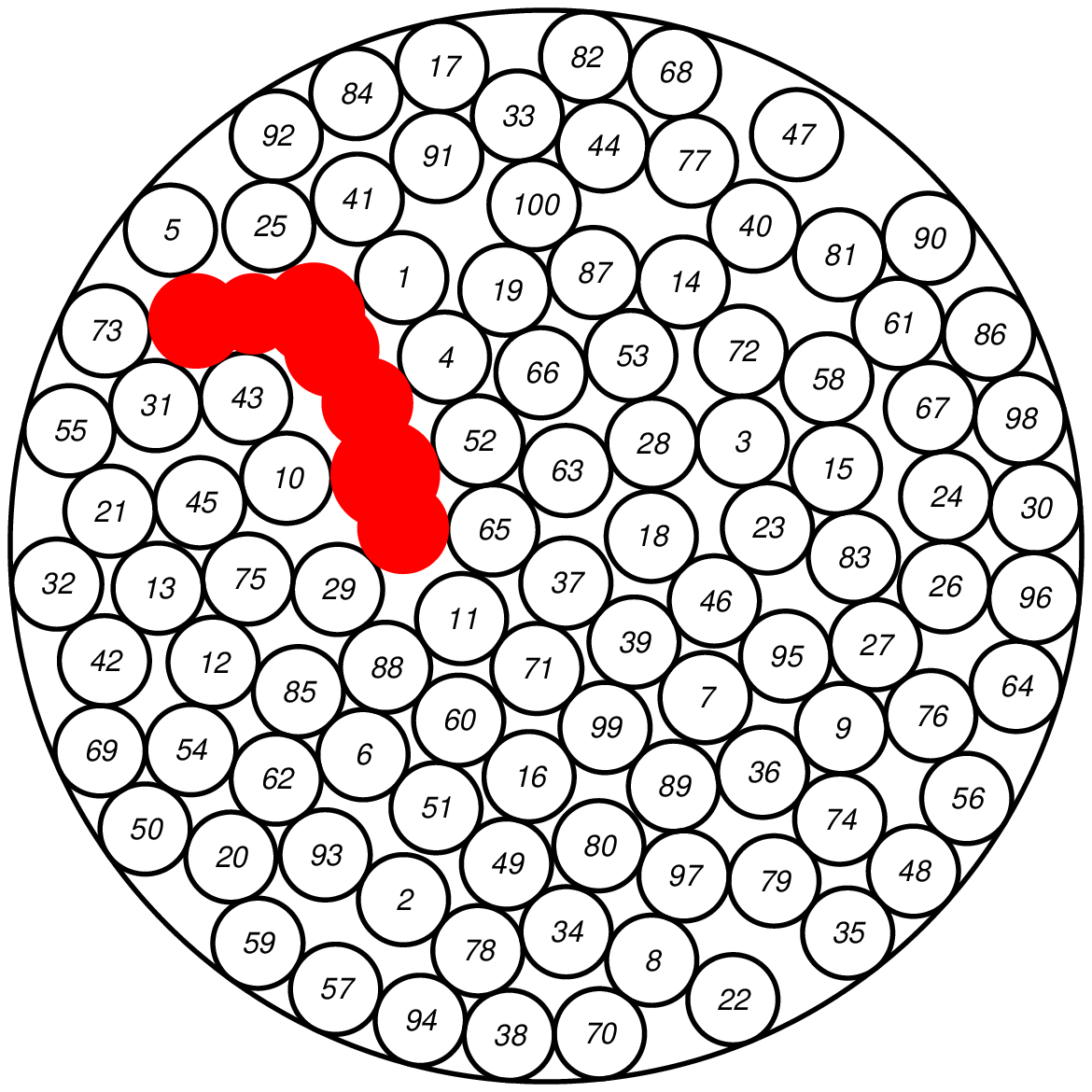}} 
\\
\subfigure[Test problem 1,  $|F|=8$]{\label{fig:100fac8}\includegraphics[scale=0.4]{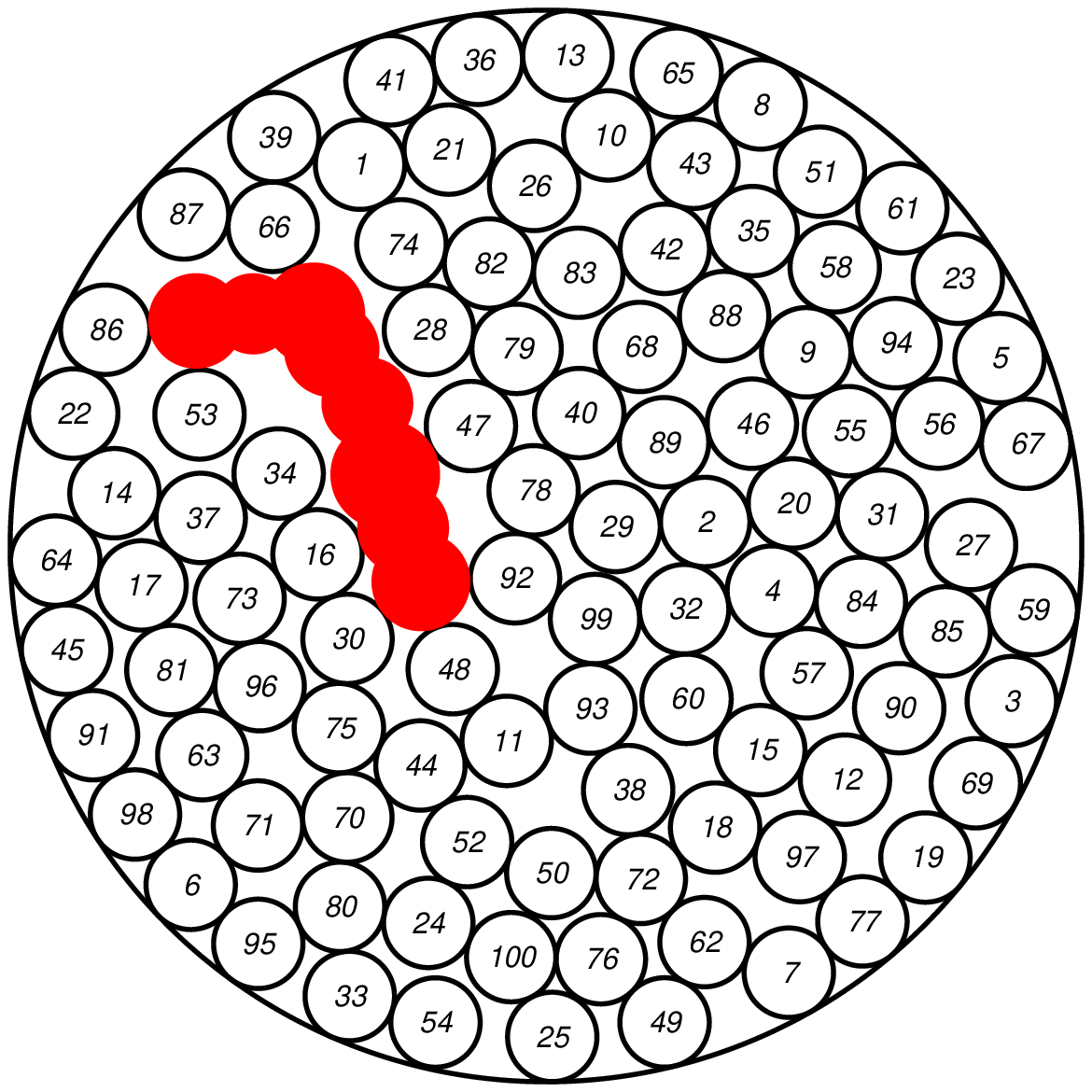}} \hspace{.2cm}
\subfigure[Test problem 1,  $|F|=9$]{\label{fig:100fac9}\includegraphics[scale=0.4]{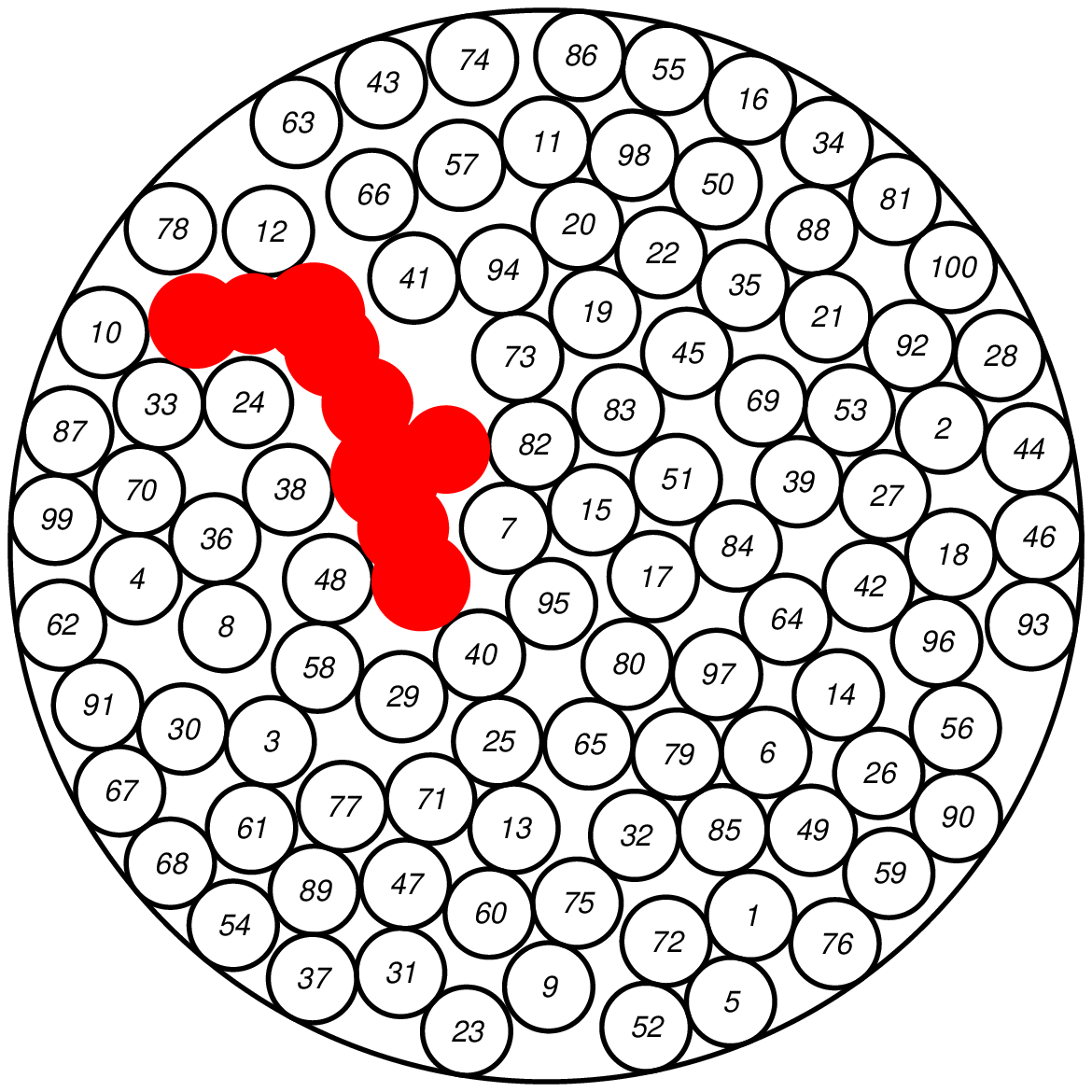}} \hspace{.2cm}
\subfigure[Test problem 1,  $|F|=10$]{\label{fig:100fac10}\includegraphics[scale=0.4]{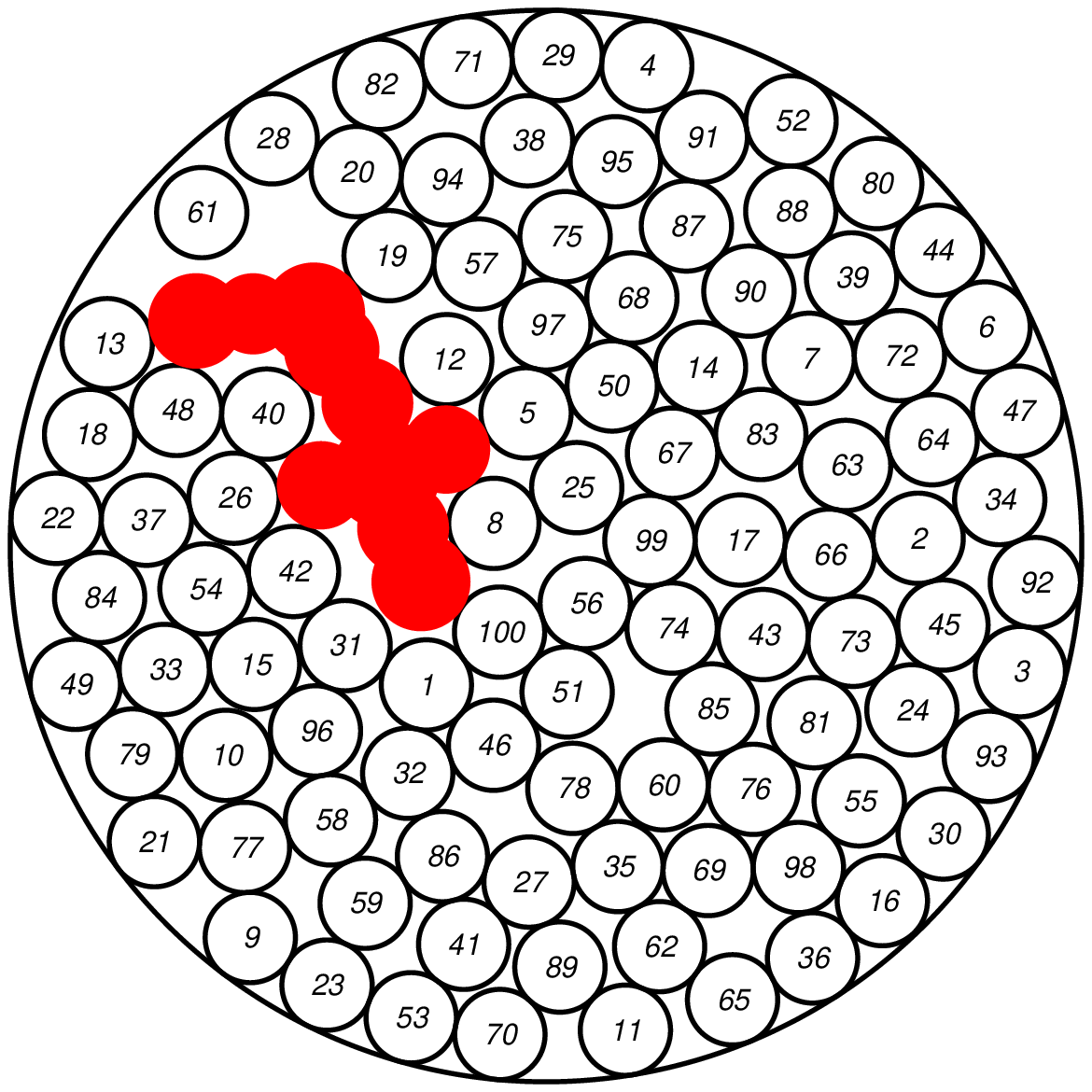}} \hspace{.2cm}
\subfigure[Test problem 1,  $|F|=11$]{\label{fig:100fa1}\includegraphics[scale=0.4]{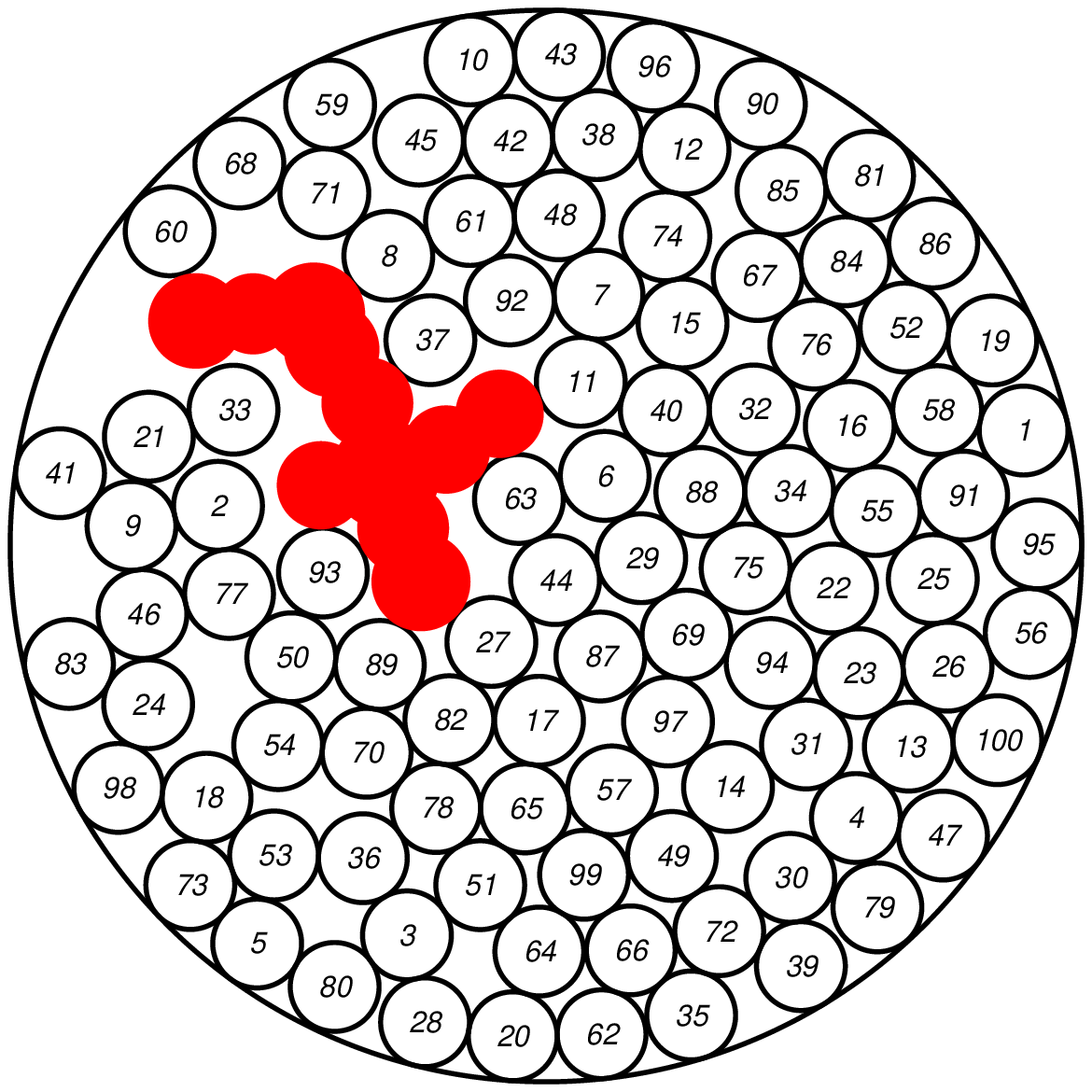}}
\label{Fig:4.100} 
\end{figure}
\end{landscape}

\begin{figure}[!htbp] 
\centering   
\caption{FSS results for test problem 6, packing 50 and 100 identical circles in a circular container with disconnected prohibited areas}
\vspace{-.5cm}
 \subfigure[$n=50$]{\label{fig:t6a}\includegraphics[scale=0.4]{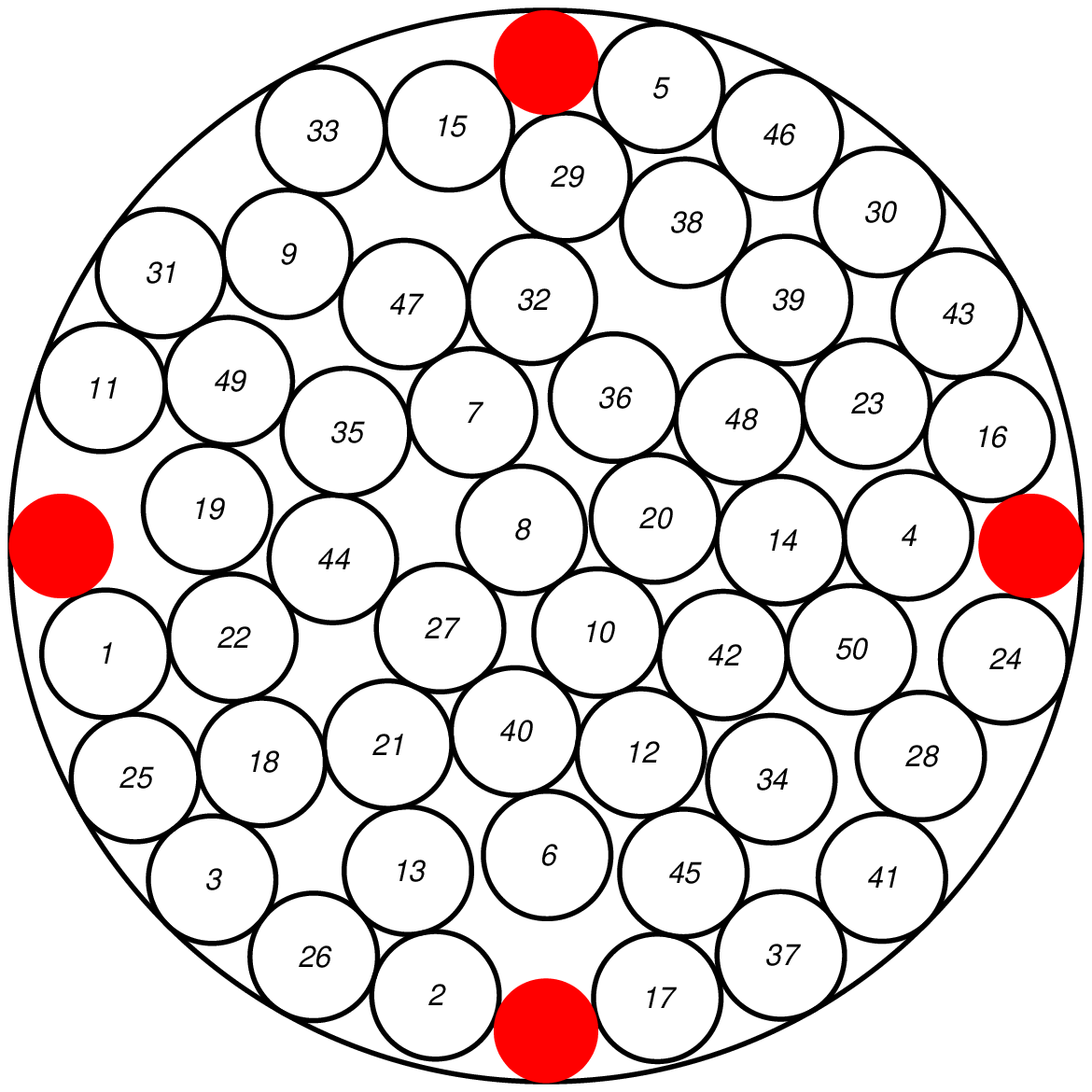}} \hspace{.2cm}
\subfigure[$n=100$]{\label{fig:t6b}\includegraphics[scale=0.4]{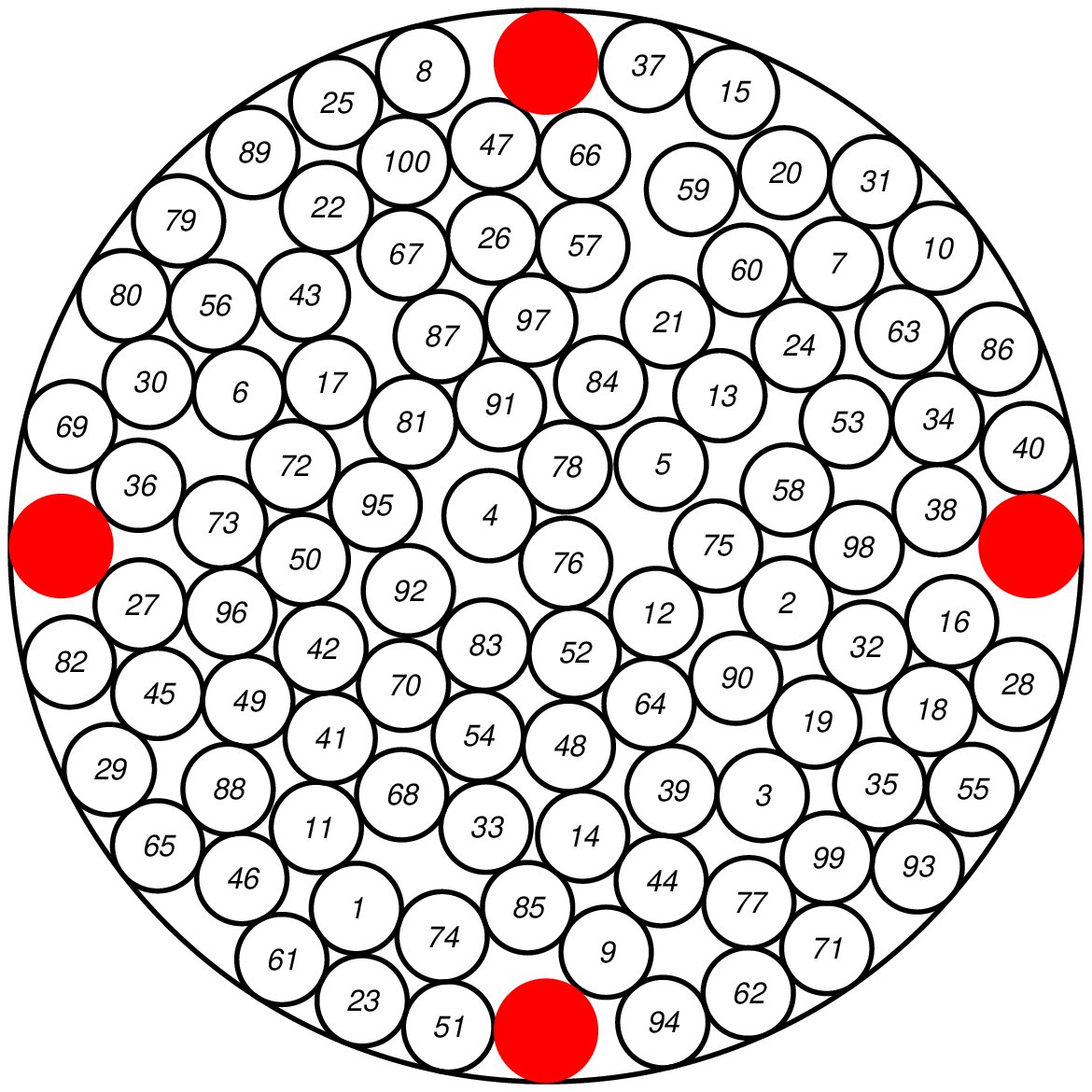}} \hspace{.2cm}
\label{Fig:test6} 
\end{figure}

\clearpage
\newpage

 \clearpage
\newpage
 \pagestyle{empty}
\linespread{1}
\small \normalsize 
\section*{\textbf{Acknowledgments}}
The first author has grant support from the programme 
UNAM-DGAPA-PAPIIT-IA106916

\end{document}